\documentclass[final,3p,times]{elsarticle}

\usepackage{amssymb}
\usepackage{amsmath}
\usepackage{amsthm}

\journal{arXiv}

\usepackage{hyperref}
\usepackage{xcolor}
\usepackage{mathtools}
\usepackage{subcaption}
\usepackage{tikz}
\usepackage{pgfplots}
\usepackage[normalem]{ulem}
\usepackage[nameinlink]{cleveref}

\hypersetup{
    colorlinks,
    linkcolor={red!50!black},
    citecolor={blue!50!black},
    urlcolor={blue!80!black}
}

\usetikzlibrary{calc}

\theoremstyle{plain}
\newtheorem{theorem}{Theorem}[section]
\newtheorem{assumption}[theorem]{Assumption}
\newtheorem{proposition}[theorem]{Proposition}
\newtheorem{remark}[theorem]{Remark}

\usepackage{etoolbox}
\AtBeginEnvironment{assumption} {\crefalias{theorem}{assumption}}
\AtBeginEnvironment{proposition}{\crefalias{theorem}{proposition}}
\AtBeginEnvironment{remark}     {\crefalias{theorem}{remark}}

\crefname{assumption}{Assumption}{Assumptions}
\Crefname{assumption}{Assumption}{Assumptions}
\crefname{proposition}{Proposition}{Propositions}
\Crefname{proposition}{Proposition}{Propositions}
\crefname{remark}{Remark}{Remarks}
\Crefname{remark}{Remark}{Remarks}

\definecolor{seaborngreen}{rgb}{0.3333333333333333, 0.6588235294117647, 0.40784313725490196}
\definecolor{seaborncyan}{rgb}{0.39215686274509803, 0.7098039215686275, 0.803921568627451}
\definecolor{seabornblue}{rgb}{0.2980392156862745, 0.4470588235294118, 0.6901960784313725}
\definecolor{seabornpurple}{rgb}{0.5058823529411764, 0.4470588235294118, 0.6980392156862745}
\definecolor{seabornred}{rgb}{0.7686274509803922, 0.3058823529411765, 0.3215686274509804}
\definecolor{seabornorange}{rgb}{0.958, 0.476, 0.206}
\definecolor{seabornsand}{rgb}{0.8, 0.7254901960784313, 0.4549019607843137}

\DeclareMathOperator{\Alt}{Alt}
\DeclareMathOperator{\sgn}{sgn}
\DeclareMathOperator{\Forall}{\forall}
\DeclareMathOperator{\tr}{tr}
\DeclareMathOperator{\sech}{sech}

\newcommand{\bbA}{\mathbb{A}}
\newcommand{\bbB}{\mathbb{B}}
\newcommand{\bbC}{\mathbb{C}}
\newcommand{\bbD}{\mathbb{D}}
\newcommand{\bbE}{\mathbb{E}}
\newcommand{\bbF}{\mathbb{F}}
\newcommand{\bbG}{\mathbb{G}}
\newcommand{\bbH}{\mathbb{H}}
\newcommand{\bbI}{\mathbb{I}}
\newcommand{\bbJ}{\mathbb{J}}
\newcommand{\bbK}{\mathbb{K}}
\newcommand{\bbL}{\mathbb{L}}
\newcommand{\bbM}{\mathbb{M}}
\newcommand{\bbN}{\mathbb{N}}
\newcommand{\bbO}{\mathbb{O}}
\newcommand{\bbP}{\mathbb{P}}
\newcommand{\bbQ}{\mathbb{Q}}
\newcommand{\bbR}{\mathbb{R}}
\newcommand{\bbS}{\mathbb{S}}
\newcommand{\bbT}{\mathbb{T}}
\newcommand{\bbU}{\mathbb{U}}
\newcommand{\bbV}{\mathbb{V}}
\newcommand{\bbW}{\mathbb{W}}
\newcommand{\bbX}{\mathbb{X}}
\newcommand{\bbY}{\mathbb{Y}}
\newcommand{\bbZ}{\mathbb{Z}}

\newcommand{\bfA}{\mathbf{A}}
\newcommand{\bfB}{\mathbf{B}}
\newcommand{\bfC}{\mathbf{C}}
\newcommand{\bfD}{\mathbf{D}}
\newcommand{\bfE}{\mathbf{E}}
\newcommand{\bfF}{\mathbf{F}}
\newcommand{\bfG}{\mathbf{G}}
\newcommand{\bfH}{\mathbf{H}}
\newcommand{\bfI}{\mathbf{I}}
\newcommand{\bfJ}{\mathbf{J}}
\newcommand{\bfK}{\mathbf{K}}
\newcommand{\bfL}{\mathbf{L}}
\newcommand{\bfM}{\mathbf{M}}
\newcommand{\bfN}{\mathbf{N}}
\newcommand{\bfO}{\mathbf{O}}
\newcommand{\bfP}{\mathbf{P}}
\newcommand{\bfQ}{\mathbf{Q}}
\newcommand{\bfR}{\mathbf{R}}
\newcommand{\bfS}{\mathbf{S}}
\newcommand{\bfT}{\mathbf{T}}
\newcommand{\bfU}{\mathbf{U}}
\newcommand{\bfV}{\mathbf{V}}
\newcommand{\bfW}{\mathbf{W}}
\newcommand{\bfX}{\mathbf{X}}
\newcommand{\bfY}{\mathbf{Y}}
\newcommand{\bfZ}{\mathbf{Z}}

\newcommand{\bfa}{\mathbf{a}}
\newcommand{\bfb}{\mathbf{b}}
\newcommand{\bfc}{\mathbf{c}}
\newcommand{\bfd}{\mathbf{d}}
\newcommand{\bfe}{\mathbf{e}}
\newcommand{\bff}{\mathbf{f}}
\newcommand{\bfg}{\mathbf{g}}
\newcommand{\bfh}{\mathbf{h}}
\newcommand{\bfi}{\mathbf{i}}
\newcommand{\bfj}{\mathbf{j}}
\newcommand{\bfk}{\mathbf{k}}
\newcommand{\bfl}{\mathbf{l}}
\newcommand{\bfm}{\mathbf{m}}
\newcommand{\bfn}{\mathbf{n}}
\newcommand{\bfo}{\mathbf{o}}
\newcommand{\bfp}{\mathbf{p}}
\newcommand{\bfq}{\mathbf{q}}
\newcommand{\bfr}{\mathbf{r}}
\newcommand{\bfs}{\mathbf{s}}
\newcommand{\bft}{\mathbf{t}}
\newcommand{\bfu}{\mathbf{u}}
\newcommand{\bfv}{\mathbf{v}}
\newcommand{\bfw}{\mathbf{w}}
\newcommand{\bfx}{\mathbf{x}}
\newcommand{\bfy}{\mathbf{y}}
\newcommand{\bfz}{\mathbf{z}}

\newcommand{\bfzero}{\mathbf{0}}
\newcommand{\bfone}{\mathbf{1}}
\newcommand{\bftwo}{\mathbf{2}}
\newcommand{\bfthree}{\mathbf{3}}
\newcommand{\bffour}{\mathbf{4}}
\newcommand{\bffive}{\mathbf{5}}
\newcommand{\bfsix}{\mathbf{6}}
\newcommand{\bfseven}{\mathbf{7}}
\newcommand{\bfeight}{\mathbf{8}}
\newcommand{\bfnine}{\mathbf{9}}

\newcommand{\bfalpha}{\bm{\alpha}}
\newcommand{\bfbeta}{\bm{\beta}}
\newcommand{\bfmu}{\bm{\mu}}
\newcommand{\bfomega}{\bm{\omega}}
\newcommand{\bfvarphi}{\bm{\varphi}}
\newcommand{\bftheta}{\bm{\theta}}

\newcommand{\rmA}{\mathrm{A}}
\newcommand{\rmB}{\mathrm{B}}
\newcommand{\rmC}{\mathrm{C}}
\newcommand{\rmD}{\mathrm{D}}
\newcommand{\rmE}{\mathrm{E}}
\newcommand{\rmF}{\mathrm{F}}
\newcommand{\rmG}{\mathrm{G}}
\newcommand{\rmH}{\mathrm{H}}
\newcommand{\rmI}{\mathrm{I}}
\newcommand{\rmJ}{\mathrm{J}}
\newcommand{\rmK}{\mathrm{K}}
\newcommand{\rmL}{\mathrm{L}}
\newcommand{\rmM}{\mathrm{M}}
\newcommand{\rmN}{\mathrm{N}}
\newcommand{\rmO}{\mathrm{O}}
\newcommand{\rmP}{\mathrm{P}}
\newcommand{\rmQ}{\mathrm{Q}}
\newcommand{\rmR}{\mathrm{R}}
\newcommand{\rmS}{\mathrm{S}}
\newcommand{\rmT}{\mathrm{T}}
\newcommand{\rmU}{\mathrm{U}}
\newcommand{\rmV}{\mathrm{V}}
\newcommand{\rmW}{\mathrm{W}}
\newcommand{\rmX}{\mathrm{X}}
\newcommand{\rmY}{\mathrm{Y}}
\newcommand{\rmZ}{\mathrm{Z}}

\newcommand{\rma}{\mathrm{a}}
\newcommand{\rmb}{\mathrm{b}}
\newcommand{\rmc}{\mathrm{c}}
\newcommand{\rmd}{\mathrm{d}}
\newcommand{\rme}{\mathrm{e}}
\newcommand{\rmf}{\mathrm{f}}
\newcommand{\rmg}{\mathrm{g}}
\newcommand{\rmh}{\mathrm{h}}
\newcommand{\rmi}{\mathrm{i}}
\newcommand{\rmj}{\mathrm{j}}
\newcommand{\rmk}{\mathrm{k}}
\newcommand{\rml}{\mathrm{l}}
\newcommand{\rmm}{\mathrm{m}}
\newcommand{\rmn}{\mathrm{n}}
\newcommand{\rmo}{\mathrm{o}}
\newcommand{\rmp}{\mathrm{p}}
\newcommand{\rmq}{\mathrm{q}}
\newcommand{\rmr}{\mathrm{r}}
\newcommand{\rms}{\mathrm{s}}
\newcommand{\rmt}{\mathrm{t}}
\newcommand{\rmu}{\mathrm{u}}
\newcommand{\rmv}{\mathrm{v}}
\newcommand{\rmw}{\mathrm{w}}
\newcommand{\rmx}{\mathrm{x}}
\newcommand{\rmy}{\mathrm{y}}
\newcommand{\rmz}{\mathrm{z}}

\newcommand{\calA}{\mathcal{A}}
\newcommand{\calB}{\mathcal{B}}
\newcommand{\calC}{\mathcal{C}}
\newcommand{\calD}{\mathcal{D}}
\newcommand{\calE}{\mathcal{E}}
\newcommand{\calF}{\mathcal{F}}
\newcommand{\calG}{\mathcal{G}}
\newcommand{\calH}{\mathcal{H}}
\newcommand{\calI}{\mathcal{I}}
\newcommand{\calJ}{\mathcal{J}}
\newcommand{\calK}{\mathcal{K}}
\newcommand{\calL}{\mathcal{L}}
\newcommand{\calM}{\mathcal{M}}
\newcommand{\calN}{\mathcal{N}}
\newcommand{\calO}{\mathcal{O}}
\newcommand{\calP}{\mathcal{P}}
\newcommand{\calQ}{\mathcal{Q}}
\newcommand{\calR}{\mathcal{R}}
\newcommand{\calS}{\mathcal{S}}
\newcommand{\calT}{\mathcal{T}}
\newcommand{\calU}{\mathcal{U}}
\newcommand{\calV}{\mathcal{V}}
\newcommand{\calW}{\mathcal{W}}
\newcommand{\calX}{\mathcal{X}}
\newcommand{\calY}{\mathcal{Y}}
\newcommand{\calZ}{\mathcal{Z}}
\newcommand{\Kn}{\mathrm{Kn}}

\usepackage{todonotes}
\newcommand{\pef}[1]{\todo[color=blue!30,inline]{{\bf PEF:} #1}}
\newcommand{\bda}[1]{\todo[color=orange!30]{{\bf BDA:} #1}}
\newcommand{\bdainline}[1]{\todo[color=orange!30,inline]{{\bf BDA:} #1}}

\newif\ifresponse
\responsefalse
\colorlet{revisedcolour}{seabornblue!50!blue}
\newcommand{\revised}[1]{\ifresponse\textcolor{revisedcolour}{#1}\else#1\fi}
\newcommand{\deleted}[1]{\ifresponse\textcolor{seabornred!50!red}{\sout{#1}}\else\fi}

\begin{document}

\begin{frontmatter}



\title{Conservative and dissipative discretisations of multi-conservative ODEs and GENERIC systems}

\author[oxford]{Boris D.~Andrews}
\author[oxford,charles]{Patrick E.~Farrell}

\affiliation[oxford]{
    organization={Mathematical Institute, University of Oxford},
    country={UK}
    }
\affiliation[charles]{
    organization={Mathematical Institute, Faculty of Mathematics and Physics, Charles University},
    country={Czechia}
    }

\begin{abstract}
Ordinary and partial differential equations describing thermodynamically isolated systems typically possess conserved quantities (like mass, momentum, and energy) and dissipated quantities (like entropy).
Preserving these conservation and dissipation laws on discretisation in time can yield vastly better approximations for the same computational effort, compared to schemes that are not structure-preserving.
In this work we present two novel contributions:
(i)~an arbitrary-order time discretisation inspired by the Nambu bracket for general conservative ordinary differential equations that conserves all \revised{prescribed} invariants, and (ii)~an energy-conserving and entropy-dissipating scheme for both ordinary and partial differential equations written in the GENERIC format, a superset of Poisson and gradient-descent systems.
In both cases the underlying strategy is the same:
the systematic introduction of auxiliary variables, allowing for the replication at the discrete level of the proofs of conservation or dissipation.
We illustrate the advantages of our approximations with numerical examples of the Kepler and Kovalevskaya problems, a combustion engine model, and the Benjamin--Bona--Mahony equation.
\end{abstract}




\end{frontmatter}


\section{Introduction}


Preserving geometric structure on discretisation has emerged as one of the key themes in modern numerical analysis.
For initial-value problems involving ordinary differential equations (ODEs) and partial differential equations (PDEs), desirable features include symplecticity, reversibility, contractivity, maximum principles, nonnegativity, conservation laws, and dissipation inequalities.
The art of designing such structure-preserving integrators is known as geometric numerical integration~\cite{Hairer_et_al_2006}.
In general, it is not possible to preserve every such feature simultaneously;
for example, the celebrated Ge--Marsden theorem shows that approximate symplectic algorithms cannot conserve energy for certain nonintegrable Hamiltonian systems~\cite{Zhong_Marsden_1988}.
To faithfully reproduce the characteristics of a problem's exact solution at the discrete level, one should target those structures that are most influential in the behaviour of the solution.

In this work, we further develop a framework recently proposed by the authors \cite{Andrews_Farrell_2025a} for devising time discretisations that exactly reproduce (up to machine precision, solver tolerances, and quadrature errors) conservation laws and dissipation inequalities.
We extend the ideas of \cite{Andrews_Farrell_2025a} in two directions.
For general conservative ODEs, we construct arbitrary-order integrators that exactly conserve all \revised{specified} invariants.
\revised{That is, given a finite collection of invariants supplied \emph{a priori}, the scheme conserves each of them exactly (up to the numerical errors noted above).}
For ODEs and PDEs written in the GENERIC (General Equation for Non-Equilibrium Reversible-Irreversible Coupling) formalism~\cite{Grmela_Ottinger_1997, Ottinger_Grmela_1997}, an abstract thermodynamically-compatible expression of certain initial-value problems, we obtain schemes that simultaneously conserve energy and dissipate entropy.
We also present numerical examples where preserving these structures appears to offer decisive advantages over typical methods (e.g.~those preserving symplecticity) for single simulations over relatively long timescales.

\revised{%
The unifying idea is the use of \emph{auxiliary variables}, i.e.~projections of certain \emph{associated test functions} (e.g.~gradients of these quantities of interest in the ODE case) onto a discrete test space, proposed in our prior work \cite{Andrews_Farrell_2025a}.
The framework of \cite{Andrews_Farrell_2025a} is general, developed abstractly and demonstrated on the incompressible and compressible Navier--Stokes equations.
The constructions and applications proposed here are novel.
This includes the alternating-form construction of \Cref{th:alternating_forms} and its use in our proposed conservative integrator \eqref{eq:integrable_avcpg}, as well as our treatment of the GENERIC formalism, including the compatibility Assumptions~\ref{ass:generic_matrices} \&~\ref{ass:generic_operators} as used in the definition of our proposed energy- and entropy-stable ODE \eqref{eq:generic_ode_avcpg} and PDE \eqref{eq:generic_pde_avcpg} integrators.
The various numerical examples are likewise all new.
}


\revised{%
For general conservative ODEs, our construction sits within a mature literature on the exact conservation of first integrals.
In the case of a single invariant, the most closely related work is the energy-preserving integrator of Cohen \& Hairer \cite{Cohen_Hairer_2011}, an arbitrary-order generalisation of the averaged-vector-field (AVF) discrete gradient method of McLachlan, Quispel \& Robidoux \cite{McLachlan_Quispel_Robidoux_1999}.
Our single-invariant scheme \eqref{eq:odecpgmod} recovers exactly this family.
Discrete multipliers, proposed by Wan, Bihlo \& Nave \cite{Wan_Bihlo_Nave_2017}, represent a further generalisation of discrete gradients, with a particular focus on non-autonomous systems, solvability, and long-term stability \cite{Wan_Nave_2018}.
}

\revised{%
For \emph{multiple} invariants, two strands of prior work are especially relevant.
Dahlby, Owren \& Yaguchi \cite{Dahlby_Owren_Yaguchi_2011} define a discrete tangent space normal to the associated discrete gradients \cite{Gonzalez_1996, McLachlan_Quispel_Robidoux_1999};
they then propose modifying an underlying Runge--Kutta method through projection onto this space, observing that this guarantees each conservation law.
In recent work, Schulz \& Wan \cite{Schulz_Wan_2025} construct the minimum-$\ell^2$-norm modification of a time discretisation that preserves orthogonality to the discrete multipliers (or discrete gradients), using a pseudoinverse.
}

\revised{%
On the line-integral side, Brugnano, Iavernaro \& Frasca-Caccia \cite{Brugnano_Iavernaro_2012, Brugnano_Iavernaro_2016, Brugnano_FrascaCaccia_Iavernaro_2019} generalise the Hamiltonian Boundary Value Method framework to conserve any number of independent invariants, augmenting the underlying polynomial method with a projection-type correction that enforces conservation of all invariants.
Each of these approaches can be shown to maintain the convergence order of the modified underlying Runge--Kutta method.
}

\revised{%
Our conservative integrator differs from the above in its construction:
we rewrite such multi-conservative systems so that their right-hand side is generated by an alternating form, contracted against the gradients of the prescribed invariants (see \Cref{th:alternating_forms}).
This leaves a natural place to introduce discrete gradients for each variable, extending to arbitrary order in a similar manner to \cite{Cohen_Hairer_2011}.
This alternating form construction bears a strong similarity to the Nambu bracket \cite{Nambu_1973, Takhtajan_1994} (see \Cref{rem:nambu}).
}

Our GENERIC discretisations sit alongside a family of works in which entropy derivatives are likewise replaced by projections onto the discrete test space.
In 2009, Romero \cite{Romero_2009} proposed introducing AVF discrete gradients for both the energy and entropy, modifying the Poisson and friction matrices to match in such a way as to preserve the GENERIC degeneracy conditions.
This guarantees energy conservation and entropy dissipation, recovering precisely our scheme \eqref{eq:generic_ode_avcpg} at lowest order in time (where these discrete gradients are interpreted as auxiliary variables).
Giesselmann, Karsai \& Tscherpel \cite{Giesselmann_Karsai_Tscherpel_2025} propose an arbitrary-order integrator for \emph{port}-Hamiltonian PDEs that preserves the dissipation of energy.
In the absence of dissipation, this reduces to an energy-conserving integrator for Hamiltonian systems;
similarly, in the absence of dissipation, our proposed GENERIC PDE integrator \eqref{eq:generic_pde_avcpg} reduces to an energy-conserving integrator for Hamiltonian systems and coincides precisely with that construction.
Our construction combines (i)~the arbitrary-order time discretisations of \cite{Giesselmann_Karsai_Tscherpel_2025} with (ii)~Poisson and friction matrices and operators modified to act on the auxiliary variables (see \Cref{ass:generic_matrices,ass:generic_operators}) to preserve the GENERIC degeneracy conditions, similar to \cite{Romero_2009}.
Lombardi \& Pagliantini \cite{Lombardi_Pagliantini_2025} consider the discretisation of Poisson systems with a gradient-descent term.
Their spatial semi-discretisation aligns closely with our proposed GENERIC integrator \eqref{eq:generic_pde_avcpg}, instead retaining the operators as evaluated on the primal variable.
For constant-valued Poisson and friction operators, their proposed semi-discretisation aligns with ours, indeed inheriting the same energy- and entropy-stability results in this case.
Moreover, the authors propose the use of AVF discrete gradients for the time discretisations, aligning with our proposed discretisation (in the constant-valued operator case) at lowest order.

A different strategy for GENERIC ODEs was taken by Shang \& \"Ottinger \cite{Shang_Ottinger_2020} who used a splitting, applying a symplectic integrator to the conservative component, and defining the friction matrix for the dissipative component in terms of the modified energy of that symplectic integrator.
This idea of defining the friction matrix in terms of a modified energy is reminiscent of our construction in \Cref{ass:generic_matrices}, albeit in the setting of symplectic integrators in place of conservative or dissipative ones.

The remainder of this manuscript proceeds as follows.
In \Cref{sec:conservative_odes}, we develop the conservative ODE integrator;
we illustrate its structure-preserving properties with numerical demonstrations on the Kepler and Kovalevskaya problems.
In \Cref{sec:generic_odes}, we present the GENERIC ODE construction, with a similar numerical example on an unfired combustion engine.
In \Cref{sec:generic_pdes}, we extend this idea to GENERIC PDEs;
to demonstrate our construction, we consider the Boltzmann and Benjamin--Bona--Mahony equations, with a numerical demonstration of the latter.
We offer some concluding remarks in \Cref{sec:conclusions}.

\section{General conservative ODEs} \label{sec:conservative_odes}

We first consider general conservative ODE systems with possibly many invariants.
Let $\bff : \bbR^d \to \bbR^d$ induce the general ODE system
\begin{equation}\label{eq:integrable_system}
    \dot{\bfx} = \bff(\bfx),
\end{equation}
with initial condition (IC) $\bfx(0) = \bfx_0$.
Suppose this system is conservative in $P \,(< d)$ independent invariants $(N_p : \bbR^d \to \bbR)_{p=1}^P$, with the property that $\nabla N_p(\bfx)^\top \bff(\bfx) = 0$ for each $p = 1, \dots, P$\revised{; here and throughout, we call a collection of invariants \emph{independent} when their gradients $(\nabla N_p)$ are linearly independent almost everywhere}.
Each $N_p$ can then be seen to be conserved over a timestep $T_n \coloneqq [t_n, t_{n+1}]$ by testing \eqref{eq:integrable_system} with $\nabla N_p$:
\begin{equation}
    N_p(x_{n+1}) - N_p(x_n) = \int_{T_n} \dot{N}_p = \int_{T_n} \nabla N_p(\bfx) \cdot \dot{\bfx} = \int_{T_n} \nabla N_p(\bfx)^\top \bff(\bfx) = 0.
\end{equation}

Before considering the general system \eqref{eq:integrable_system} in further detail, it is instructive to consider a Poisson system,
\revised{%
\begin{subequations} \label{eq:poisson}
\begin{equation}\label{eq:poisson_bracket}
    \dot{\bfx} = \bff(\bfx) \coloneqq \{\bfx, H\},
\end{equation}
with Poisson bracket $\{\cdot,\cdot\}$ and conserved Hamiltonian $H : \bbR^d \to \bbR$.
More amenable to discretisation (in particular, to the later introduction of a discrete gradient approximating $\nabla H$) is the equivalent \emph{matrix formulation}
\begin{equation}\label{eq:poisson_tensor}
    \dot{\bfx} = \bff(\bfx) \coloneqq B(\bfx) \nabla H(\bfx),
\end{equation}
\end{subequations}
}
where $B(\bfx) \in \bbR^{d \times d}$ is a skew-symmetric matrix encoding the Poisson bracket\deleted{, and $H$ is the conserved Hamiltonian}\footnote{
    Since we do not impose the Jacobi identity on $B$, this would more precisely be referred to as \emph{almost} Poisson.
    However, we use the former terminology throughout this manuscript for brevity.
}.
In this case, the relation $\nabla H(\bfx)^\top \bff(\bfx) = \nabla H(\bfx)^\top B(\bfx) \nabla H(\bfx) = 0$ follows from the skew-symmetry of $B$.
A broad class of discretisations for \eqref{eq:poisson} is of the form:
find $\bfx \in \bbX_n$ such that
\begin{equation} \label{eq:conservative:base_scheme}
    \calI_n\left[ \bfy^\top \dot{\bfx} \right] = \calI_n\left[\bfy^\top B(\bfx) \nabla H(\bfx)\right]
\end{equation}
for all $\bfy \in \dot{\bbX}_n$. Here $\bbX_n$ is the trial space of polynomials of a fixed degree $s \ge 1$ satisfying known initial data at the beginning of the timestep,
\begin{equation} \label{eq:solution_space_ode}
    \bbX_n \coloneqq \left\{ \, \bfy \in \bbP_s(T_n; \bbR^d) \mid \bfy(t_n) \text{ satisfies known initial data} \, \right\},
\end{equation}
with its time derivative satisfying $\dot{\bbX}_n = \bbP_{s-1}(T_n; \bbR^d)$.
The operator $\calI_n$ is a linear approximation to the integral over $T_n$ (a quadrature rule), i.e.~$\calI_n[\phi] \approx \int_{T_n} \phi$, specifying the time discretisation;
the approximation must be sign-preserving, i.e.
\begin{subequations}\label{eq:In}
\begin{equation}\label{eq:Isignpreserving}
    \phi \ge 0 \implies \calI_n[\phi] \ge 0,
\end{equation}
appropriately scaled in $\Delta t_n \coloneqq t_{n+1} - t_n$, i.e.
\begin{equation}\label{eq:Iscaling}
    \calI_n[1] = \Delta t_n,
\end{equation}
\end{subequations}
and the map $\phi \mapsto \calI_n[\phi^2]^\frac{1}{2}$ must define a norm on $\bbP_{s-1}(T_n)$.
Examples of such linear operators include any $\ge s$-stage quadrature rule with positive weights, and the exact integral.
In the case where $\calI_n$ is an $s$-stage quadrature rule, \eqref{eq:conservative:base_scheme} is equivalent to a collocation method at the quadrature points (see \cite[Sec.~II.1.2]{Hairer_Lubich_Wanner_2006});
in the case where $\calI_n$ is the exact integral, \eqref{eq:conservative:base_scheme} is a continuous Petrov--Galerkin discretisation in time (see \cite[Chap.~70]{Ern_Guermond_2021c}).

We wish to discretise \eqref{eq:poisson} in time while retaining the conservation of $H$.
However, the time discretisation \eqref{eq:conservative:base_scheme} cannot in general achieve this.
The central difficulty with reproducing the conservation argument for the continuous case is that $\nabla H(\bfx) \not\in \dot{\bbX}_n$, i.e.~$\nabla H(\bfx)$ is not a valid choice of discrete test function.
This may be overcome by introducing an auxiliary variable $\widetilde{\nabla H} \in \dot{\bbX}_n$ approximating $\nabla H(\bfx)$, yielding the scheme:
find $(\bfx, \widetilde{\nabla H}) \in \bbX_n \times \dot{\bbX}_n$ such that
\begin{subequations}\label{eq:odecpgmod}
\begin{align}
    \calI_n\left[\bfy^\top \dot{\bfx}\right]  &=  \calI_n\left[\bfy^\top B(\bfx) \widetilde{\nabla H}\right], \label{eq:odecpgmoda} \\
    \calI_n\left[\widetilde{\nabla H}^\top\bfy_0\right]  &=  \int_{T_n}{\nabla H}(\bfx)^\top\bfy_0, \label{eq:odecpgmodb}
\end{align}
\end{subequations}
for all $(\bfy, \bfy_0) \in \dot{\bbX}_n \times \dot{\bbX}_n$.

\begin{proposition}[Conservation properties of \eqref{eq:odecpgmod}]\label{th:H_conservation}
    Assuming solutions exist, the integrator \eqref{eq:odecpgmod} satisfies $H(\bfx(t_{n+1})) = H(\bfx(t_n))$ up to quadrature errors, solver tolerances, and machine precision \revised{(see \Cref{rem:propagation})} at every timestep.
\end{proposition}

\begin{proof}
    Mimicking the continuous argument, we consider respectively $\bfy_0 = \dot{\bfx}$ in \eqref{eq:odecpgmodb} and $\bfy = \widetilde{\nabla H}$ in \eqref{eq:odecpgmoda}:
    \begin{equation}
        H(\bfx(t_{n+1})) - H(\bfx(t_n))
        =  \int_{T_n} \dot{H}
        =  \int_{T_n}\nabla H(\bfx)^\top\dot{\bfx}
        =  \calI_n\left[\widetilde{\nabla H}^\top\dot{\bfx}\right]
        =  \calI_n\left[\widetilde{\nabla H}^\top B(\bfx) \widetilde{\nabla H}\right]
        =  0.
    \end{equation}
\end{proof}

Summarising the strategy in the construction of \eqref{eq:odecpgmod}, we (i)~introduce an \emph{auxiliary variable} corresponding to the quantity we wish to preserve, namely a projection of its gradient onto the discrete test space $\dot{\bbX}_n$, and (ii)~modify the right-hand side to use the auxiliary variable so that the conservation proof can be replicated discretely.
In the context of \cite{Andrews_Farrell_2025a}, the gradient $\nabla H$ is the \emph{associated test function} corresponding to $H$;
namely, it is the function against which one would test the equation \eqref{eq:poisson} to derive the conservation or dissipation law in the continuous setting, and consequently that which we introduce into our discretisation through a discrete auxiliary variable to preserve that law.

\begin{remark}\label{rem:cohen_hairer}
    In practice, the scheme \eqref{eq:odecpgmod} incurs little additional computational cost over \eqref{eq:conservative:base_scheme} as the auxiliary variable $\widetilde{\nabla H}$ can be evaluated explicitly.
    In the case of $\calI_n$ an $s$-stage quadrature rule, the resulting reformulation can then be shown to be equivalent to the scheme of Cohen \& Hairer \cite{Cohen_Hairer_2011}, which is in turn at lowest order in time equivalent to the AVF discrete-gradient method of McLachlan, Quispel \& Robidoux \cite{McLachlan_Quispel_Robidoux_1999}.
    However, our interpretation extends more easily to the conservation of additional invariants, alongside the consideration of PDEs wherein these auxiliary variables cannot in general be evaluated explicitly.
\end{remark}

\revised{%
When $\calI_n = \int_{T_n}$ (i.e.~exact quadrature is used) the scheme \eqref{eq:odecpgmod} reduces to a continuous Petrov--Galerkin discretisation\footnote{
    \revised{%
    It is not immediately apparent that this is a continuous Petrov--Galerkin problem in the classical sense, due to the auxiliary variable $\widetilde{\nabla H}$ being a discontinuous-in-time variable of one lower degree.
    The simplest solution is to suppose one is solving equivalently for $\int_{t_n}^t \widetilde{\nabla H}(s) \rmd s \in \bbX$, substituting $\widetilde{\nabla H} = \frac{\rmd}{\rmd t}[\int_{t_n}^t \widetilde{\nabla H}(s) \rmd s]$.
    This removes the algebraic nature of $\widetilde{\nabla H}$ by reparametrisation.
    Thus classical results for Galerkin-in-time discretisations are applicable.
    }
} of a mixed formulation of the Poisson problem \eqref{eq:poisson}.
The analysis is therefore largely classical, including the existence of solutions for sufficiently small timesteps, and $\calO(\Delta t^{2s})$ convergence at mesh points \cite{Hulme_1972a}, similar to a Gau\ss{} method.
When $\calI_n \ne \int_{T_n}$ the convergence order is limited by the approximation order $r$ of $\calI_n$ to $\calO(\Delta t^{\min\{r+1,2s\}})$ \cite{Hulme_1972b}.
Such convergence results were shown by Cohen \& Hairer when $\calI_n$ is an $s$-stage quadrature rule (see \Cref{rem:cohen_hairer}) by interpreting the scheme as a partitioned Runge--Kutta method \cite{Cohen_Hairer_2011};
in this case the convergence order $\calO(\Delta t^{r+1})$ corresponds with that of the associated collocation Runge--Kutta method.
By the same methodology, similar existence and convergence results are to be expected (and observed in \Cref{sec:kepler_convergence}) for the multi-conservative \eqref{eq:integrable_avcpg} and GENERIC \eqref{eq:generic_ode_avcpg} ODE integrators proposed below.
In each case appropriate regularity assumptions are required \cite[Sec.~6.1.1]{Andrews_2025}.
}

\begin{remark}\label{rem:propagation}
    \revised{%
    For exact solutions of \eqref{eq:odecpgmod}, this discrete conservation (\Cref{th:H_conservation}) holds exactly.
    The recurring qualifier ``up to quadrature errors, solver tolerances, and machine precision'' reflects that \eqref{eq:odecpgmod}, and our other proposed schemes below, cannot typically be solved exactly in practice, due to a combination of inexact quadrature for the exact integrals $\int_{T_n}$, finite solver tolerances, and finite-precision arithmetic.
    As these per-step errors are not corrected, they may (as for any such method) accumulate over the course of the simulation, at worst linearly in the number of timesteps.
    For the numerical examples in this manuscript, the propagated conservation errors remain negligible over the considered timescales, especially in comparison to the errors generated by non-conserving methods.
    A simple remedy to this (at the cost of symmetry in time, and not implemented in this work) would be to re-project the discrete solution onto the level sets of the invariants at each timestep \cite[Sec.~IV.4]{Hairer_Lubich_Wanner_2006}.
    }
\end{remark}

We now return to the general case \eqref{eq:integrable_system}. The same strategy suggests (i)~introducing auxiliary variables $(\widetilde{\nabla N}_p)$ approximating each $\nabla N_p(\bfx)$, and (ii)~using these new variables in the right-hand side of \eqref{eq:integrable_system}, as done in \eqref{eq:odecpgmod}.
However, as written, \eqref{eq:integrable_system} does not appear to depend on $(\nabla N_p)$.
The following theorem demonstrates that in fact $\bff$ may be rewritten to expose an explicit dependence on $(\nabla N_p)$, thereby enabling the incorporation of the auxiliary variables $(\widetilde{\nabla N}_p)$.

This result is fully constructive.
In particular, we rewrite $\bff$ in terms of an \emph{alternating form}, i.e.~a multilinear map $F : V^n \to \bbR$ over a vector space $V$ such that $F[v_1, \dots, v_n] = 0$ whenever $v_i = v_j$ for some $i \ne j$.
We denote the space of alternating $n$-forms over $V$ by $\Alt^n V$, and define the alternatisation $\Alt F \in \Alt^n V$ of an $n$-multilinear map by
\begin{equation}\label{eq:alternatisation}
    \Alt F[v_1, \dots, v_n]  \coloneqq  \sum_{\sigma \in S_n}\sgn_\sigma F[v_{\sigma_1}, \dots, v_{\sigma_n}],
\end{equation}
where $S_n$ denotes the permutation group of degree $n$, and $\sgn_\sigma \in \{\pm 1\}$ the sign of $\sigma \in S_n$ \cite{Tu_2010}.

\begin{theorem}[Identification of alternating forms]\label{th:alternating_forms}
    For the general conservative system \eqref{eq:integrable_system} there exists $\tilde{F} : \bbR^d \to \Alt^{P+1}\bbR^d$ such that $\Forall \bfx, \bfy \in \bbR^d$,
    \begin{equation}\label{eq:alternating_form_coincidence}
        \bfy^\top\bff(\bfx)  =  \tilde{F}(\bfx)[\nabla N_1(\bfx), \dots, \nabla N_P(\bfx), \bfy].
    \end{equation}
\end{theorem}

\begin{proof}
    We demonstrate the existence of $\tilde{F}$ by construction.
    Through the independence of $(N_p)$, the gradients $(\nabla N_p)$ are linearly independent almost everywhere.
    Consequently, we may define a dual basis $(\bfm_q : \bbR^d \to \bbR^d)_{q=1}^P$ such that almost everywhere $\nabla N_p(\bfx)^\top\bfm_q(\bfx) = \delta_{pq}$.
    For each $\bfx \in \bbR^d$, define the multilinear map $\tilde{G}(\bfx) : (\bbR^d)^{P+1} \to \bbR$,
    \begin{equation}
        \tilde{G}(\bfx)[\bfn_1, \dots, \bfn_P, \bfy]  \coloneqq  \left(\bfn_1^\top\bfm_1(\bfx)\right)\cdots\left(\bfn_P^\top\bfm_P(\bfx)\right)\left(\bfy^\top\bff(\bfx)\right).
    \end{equation}
    By the orthogonality property $\nabla N_p(\bfx)^\top\bfm_q(\bfx) = \delta_{pq}$, we observe that
    \begin{equation}
        \tilde{G}(\bfx)[\nabla N_1(\bfx), \dots, \nabla N_P(\bfx), \bfy]  =  \bfy^\top\bff(\bfx);
    \end{equation}
    furthermore, by the orthogonality $\nabla N_p(\bfx)^\top\bff(\bfx) = 0$ (inherent in the conservation of $N_p$) this evaluates to zero under any non-trivial permutation of the arguments.
    Now define $\tilde{F}(\bfx) \coloneqq \Alt\tilde{G}(\bfx) \in \Alt^{P+1}\bbR^d$ to be the alternatisation of $\tilde{G}(\bfx)$.
    This $\tilde{F}(\bfx)$ is alternating for all arguments by construction, and coincides with $\tilde{G}(\bfx)$ when evaluated at $[\nabla N_1(\bfx), \dots, \nabla N_P(\bfx), \bfy]$ for any $\bfy$, since all but the trivial permutation evaluate to zero in the alternatisation~\eqref{eq:alternatisation}.
    Hence \eqref{eq:alternating_form_coincidence} holds.
\end{proof}

Since the proof here is constructive, one may potentially use it directly when seeking to define such an $\tilde{F}$.
In simpler cases, however, such an $\tilde{F}$ can often be found by inspection, as we will demonstrate for the Kepler problem in \Cref{sec:kepler} below.

\begin{remark}\label{rem:nambu}
    \revised{%
    \emph{Nambu brackets} \cite{Nambu_1973, Takhtajan_1994} represent a generalisation of the Poisson bracket, generating dynamics from several invariants $N_1, \dots, N_P$:
    \begin{equation}\label{eq:nambu_bracket}
        \dot{\bfx} = \bff(\bfx) \coloneqq \{\bfx, N_1, \dots, N_P\}.
    \end{equation}
    Like the Poisson bracket, the Nambu bracket $\{\cdot, \dots, \cdot\}$ is a multilinear and antisymmetric function of the gradients of its arguments.
    Just as the matrix formulation \eqref{eq:poisson_tensor} of the Poisson system relates to the bracket formulation \eqref{eq:poisson_bracket}, the alternating form formulation \eqref{eq:alternating_form_coincidence} of the multi-conservative system relates to the Nambu bracket formulation \eqref{eq:nambu_bracket}.
    In the same way that the skew-symmetric matrix $B(\bfx) \in \bbR^{d\times d}$ encodes the Poisson bracket, the alternating form $\tilde{F}(\bfx) \in \Alt^{P+1} \bbR^d$ encodes the Nambu bracket.
    }
\end{remark}

With \Cref{th:alternating_forms} established and $\tilde{F}$ defined, we may now apply the same general strategy described in \cite{Andrews_Farrell_2025a} to construct an integrator for \eqref{eq:integrable_system} that preserves all \revised{$P$ prescribed} conservation laws.
We introduce auxiliary variables $(\widetilde{\nabla N}_p) \in \dot{\bbX}_n^P$, approximating $(\nabla N_p(\bfx))$, such that
\begin{equation}
    \calI_n\left[\widetilde{\nabla N}_p^\top \bfy_p\right]  =  \int_{T_n} \nabla N_p(\bfx)^\top \bfy_p,
\end{equation}
for all $(\bfy_p)_{p=1}^P \in \dot{\bbX}_n^P$. The general conservative scheme is then:
find $(\bfx, (\widetilde{\nabla N}_p)) \in \bbX_n \times \dot{\bbX}_n^P$ such that
\begin{subequations}\label{eq:integrable_avcpg}
\begin{align}
    \calI_n\left[\bfy^\top \dot{\bfx}\right]               &=  \calI_n\left[\tilde{F}(\bfx)[\widetilde{\nabla N}_1, \dots, \widetilde{\nabla N}_P, \bfy]\right], \label{eq:integrable_avcpg_a}\\
    \calI_n\left[\widetilde{\nabla N}_p^\top\bfy_p\right]  &=  \int_{T_n}{\nabla N}_p(\bfx)^\top\bfy_p,  \qquad  p = 1, \dots, P,  \label{eq:integrable_avcpg_b}
\end{align}
\end{subequations}
for all $(\bfy, (\bfy_p)) \in \dot{\bbX}_n \times \dot{\bbX}_n^P$.

\begin{theorem}[Conservation properties of \eqref{eq:integrable_avcpg}]\label{th:integrable_sp}
    Assuming solutions exist, the integrator \eqref{eq:integrable_avcpg} satisfies $N_p(\bfx(t_{n+1})) = N_p(\bfx(t_n))$ for all $p$.
\end{theorem}

\begin{proof}
    For each $p$, by considering respectively $\bfy_p = \dot{\bfx}$ in \eqref{eq:integrable_avcpg_b} and $\bfy = \widetilde{\nabla N}_p$ in \eqref{eq:integrable_avcpg_a},
    \begin{equation}
        N_p(\bfx(t_{n+1})) - N_p(\bfx(t_n))
            = \int_{T_n} \dot{N}_p
            = \int_{T_n} \nabla N_p(\bfx)^\top \dot{\bfx}
            = \calI_n\left[\widetilde{\nabla N}_p^\top \dot{\bfx}\right]
            = \calI_n\left[\tilde{F}(\bfx)[\widetilde{\nabla N}_1, \dots, \widetilde{\nabla N}_P, \widetilde{\nabla N}_p]\right]
            = 0,
    \end{equation}
    where the final equality holds by the alternating property of $\tilde{F}(\bfx)$.
\end{proof}

\begin{remark}
As with \eqref{eq:odecpgmod}, the auxiliary variables introduced in \eqref{eq:integrable_avcpg} can always be eliminated, and thus incur very little additional computational cost.
\end{remark}

\subsection{Kepler problem}\label{sec:kepler}

As a numerical demonstration of the scheme \eqref{eq:integrable_avcpg} we discretise the nondimensionalised two-body Kepler problem,
\begin{align}\label{eq:kepler}
    \dot{\bfx}  =  \bfv,  \quad \quad
    \dot{\bfv}  =  - \frac{1}{\|\bfx\|^3}\bfx,
\end{align}
for $\bfx, \bfv : \bbR_+ \to \bbR^d$ representing the position and velocity respectively, and $\|\cdot\|$ denoting the $\ell^2$ norm.
Trajectories of \eqref{eq:kepler} preserve the energy $H$, angular momentum $\bfL$, and Runge--Lenz vector $\bfA$, defined as
\begin{equation}
    H(\bfx, \bfv)  \coloneqq  \frac{1}{2}\|\bfv\|^2 - \frac{1}{\|\bfx\|},  \qquad
    \bfL(\bfx, \bfv)  \coloneqq  \bfx \times \bfv,  \qquad
    \bfA(\bfx, \bfv)  \coloneqq  \bfv \times \bfL(\bfx, \bfv) - \frac{1}{\|\bfx\|}\bfx,
\end{equation}
where $\times$ denotes the cross product.
Roughly speaking, $H$ and $\bfL$ encode the shape of the orbit and the plane to which it is restricted, whereas the orientation of the orbit within that plane is encoded in $\bfA$ \cite{Taff_1985}.
These invariants are not independent, as $\|\bfA\|^2 = 1 + 2H\|\bfL\|^2$, while $\bfA$ and $\bfL$ are necessarily orthogonal;
these 3 invariants thus represent $2d-1$ independent constants of motion (for $d \in \{2, 3\}$)\deleted{, the maximum possible number of conserved quantities}.
\revised{%
A $d$-dimensional Hamiltonian system admits at most $2d-1$ independent constants of motion, one fewer than the dimension $2d$ of the phase space;
systems attaining this bound (such as the Kepler problem) are termed \emph{maximally superintegrable} \cite{Evans_1990, Miller_Post_Winternitz_2013}.}

We consider the two-dimensional case $d = 2$.
In this case, if $H$ and $\bfA = (A_1, A_2)$ are conserved, then the scalar angular momentum $L = \|\bfL\|$  will be conserved automatically, since $\|\bfA\|^2 = 1 + 2HL^2$.
We may therefore construct a fully conservative numerical integrator for $d = 2$ using our scheme \eqref{eq:integrable_avcpg} by conserving the $P=3$ independent invariants $H$, $A_1$, and $A_2$.

To apply \eqref{eq:integrable_avcpg} we must construct some $\tilde{F} : \bbR^{2 \times 2} \to \Alt^4 \bbR^{2 \times 2}$ satisfying the conditions of \eqref{eq:alternating_form_coincidence}, i.e.~such that $\Forall (\bfx, \bfv), (\bfy, \bfw) \in \bbR^{2 \times 2}$,
\begin{equation}
    \bfy^\top \bfv - \bfw^\top \frac{1}{\|\bfx\|^3}\bfx
    =
    \tilde{F}\left(\begin{pmatrix}
        \bfx \\ \bfv
    \end{pmatrix}\right)\left[\begin{pmatrix}
        \nabla_\bfx H \\ \nabla_\bfv H
    \end{pmatrix}, \begin{pmatrix}
        \nabla_\bfx A_1 \\ \nabla_\bfv A_1
    \end{pmatrix}, \begin{pmatrix}
        \nabla_\bfx A_2 \\ \nabla_\bfv A_2
    \end{pmatrix}, \begin{pmatrix}
        \bfy \\ \bfw
    \end{pmatrix}\right],
\end{equation}
where $\nabla_\bfx$, $\nabla_\bfv$ denote gradients with respect to $\bfx$, $\bfv$ respectively.
Instead of using the constructive proof in~\Cref{th:alternating_forms}, we may more simply note that the space $\Alt^4 \bbR^4$ only has dimension $1$;
any alternating $n$-form in $n$ dimensions is some multiple of the determinant map on the $n$-by-$n$ square matrix formed by the $n$ argument vectors.
This allows us to substantially reduce the space of potential maps $\tilde{F}$ to consider.
Noting the gradients of our conserved quantities,
\begin{subequations}
\begin{align}
    \nabla_\bfx H    &=  \frac{1}{\|\bfx\|^3}\bfx,  &
    \nabla_\bfx\bfA  &=  \frac{1}{\|\bfx\|^3}\bfx\otimes\bfx - \bfv\otimes\bfv + \left(\|\bfv\|^2 - \frac{1}{\|\bfx\|}\right) I,  \\
    \nabla_\bfv H    &=  \bfv,  &
    \nabla_\bfv\bfA  &=  2\bfx\otimes\bfv - \bfv\otimes\bfx - (\bfx\cdot\bfv) I,
\end{align}
\end{subequations}
where $\otimes$ denotes the outer product, we may see by inspection that, for all $(\bfx, \bfv), (\bfy, \bfw) \in \bbR^{2 \times 2}$,
\begin{equation}\label{eq:kepler_alternating_form_coincidence}
    \frac{1}{2L(\bfx, \bfv)H(\bfx, \bfv)}\det\left[\begin{matrix}
        \bfy  &  \nabla_\bfx H  &  \nabla_\bfx\bfA^\top  \\
        \bfw  &  \nabla_\bfv H  &  \nabla_\bfv\bfA^\top
    \end{matrix}\right]
    =
    \bfy^\top \bfv - \bfw^\top \frac{1}{\|\bfx\|^3}\bfx,
\end{equation}
where $\det : \bbR^{4 \times 4} \to \bbR$ denotes the determinant map.
We may therefore define $\tilde{F} : \bbR^{2 \times 2} \to \Alt^4\bbR^{2 \times 2}$ as
\begin{equation}\label{eq:kepler_alternating_form}
    \tilde{F}\left(\begin{pmatrix}
        \bfx \\ \bfv
    \end{pmatrix}\right)\left[\begin{pmatrix}
        \widetilde{\nabla_\bfx H} \\ \widetilde{\nabla_\bfv H}
    \end{pmatrix}, \begin{pmatrix}
        \widetilde{\nabla_\bfx A}_1 \\ \widetilde{\nabla_\bfv A}_1
    \end{pmatrix}, \begin{pmatrix}
        \widetilde{\nabla_\bfx A}_2 \\ \widetilde{\nabla_\bfv A}_2
    \end{pmatrix}, \begin{pmatrix}
        \bfy \\ \bfw
    \end{pmatrix}\right]
    \coloneqq
    \frac{1}{2L(\bfx, \bfv)H(\bfx, \bfv)}\det\left[\begin{matrix}
        \bfy  &  \widetilde{\nabla_\bfx H}  &  \widetilde{\nabla_\bfx A}_1  &  \widetilde{\nabla_\bfx A}_2  \\
        \bfw  &  \widetilde{\nabla_\bfv H}  &  \widetilde{\nabla_\bfv A}_1  &  \widetilde{\nabla_\bfv A}_2
    \end{matrix}\right],
\end{equation}
which satisfies \eqref{eq:alternating_form_coincidence} by \eqref{eq:kepler_alternating_form_coincidence}, and can be seen to be an alternating form by the alternating properties of the determinant.

Through \eqref{eq:integrable_avcpg} we then arrive at our fully conservative integrator for the two-dimensional Kepler problem:
find $((\bfx, \bfv), (\widetilde{\nabla_\bfx H}, \widetilde{\nabla_\bfv H}), (\widetilde{\nabla_\bfx \bfA}, \widetilde{\nabla_\bfv \bfA})) \in \bbX_n \times \dot{\bbX}_n \times \dot{\bbX}_n^2$ (with $\bbX_n$ defined as in \eqref{eq:solution_space_ode} for $d = 2 \times 2$) such that
\begin{subequations}\label{eq:kepler_avcpg}
\begin{align}
    \calI_n\left[\bfy^\top \dot{\bfx} + \bfw^\top \dot{\bfv}\right]
        &=  \calI_n\left[\frac{1}{2L(\bfx, \bfv)H(\bfx, \bfv)}\det\left[\begin{matrix}
            \bfy  &  \widetilde{\nabla_\bfx H}  &  \widetilde{\nabla_\bfx \bfA}^\top  \\
            \bfw  &  \widetilde{\nabla_\bfv H}  &  \widetilde{\nabla_\bfv \bfA}^\top
        \end{matrix}\right]\right],  \\
    \calI_n\left[\widetilde{\nabla_\bfx H}^\top \bfy_H + \widetilde{\nabla_\bfv H}^\top \bfw_H\right]
        &=  \int_{T_n} \nabla H(\bfx, \bfv)^\top \bfy_H + \nabla H(\bfx, \bfv)^\top \bfw_H,  \\
    \calI_n\left[\tr\Big(\widetilde{\nabla_\bfx \bfA} Y_A + \widetilde{\nabla_\bfv \bfA} W_A\Big)\right]
        &=  \int_{T_n} \tr\Big(\nabla_\bfx \bfA(\bfx, \bfv) Y_A + \nabla_\bfv \bfA(\bfx, \bfv) W_A\Big),
\end{align}
\end{subequations}
for all $((\bfy, \bfw), (\bfy_H, \bfw_H), (Y_A, W_A)) \in \dot{\bbX}_n \times \dot{\bbX}_n \times \dot{\bbX}_n^2$, where $\tr$ denotes the trace.

We test the fully conservative integrator \eqref{eq:kepler_avcpg} numerically on a standard set of ICs (inspired by \cite[Sec.~I.2.3]{Hairer_Lubich_Wanner_2006}) given by $\bfx(0) = (0.4, 0)$, $\bfv(0) = (0, 2)$.
In each of these tests, we take $\calI_n$ to be an $s$-stage Gau\ss{}--Legendre quadrature, such that the scheme \eqref{eq:kepler_avcpg} is a conservative modification of a Gau\ss{} method of the same order.

\subsubsection{Comparison test}\label{sec:kepler_comparison}

\begin{figure}[!ht]
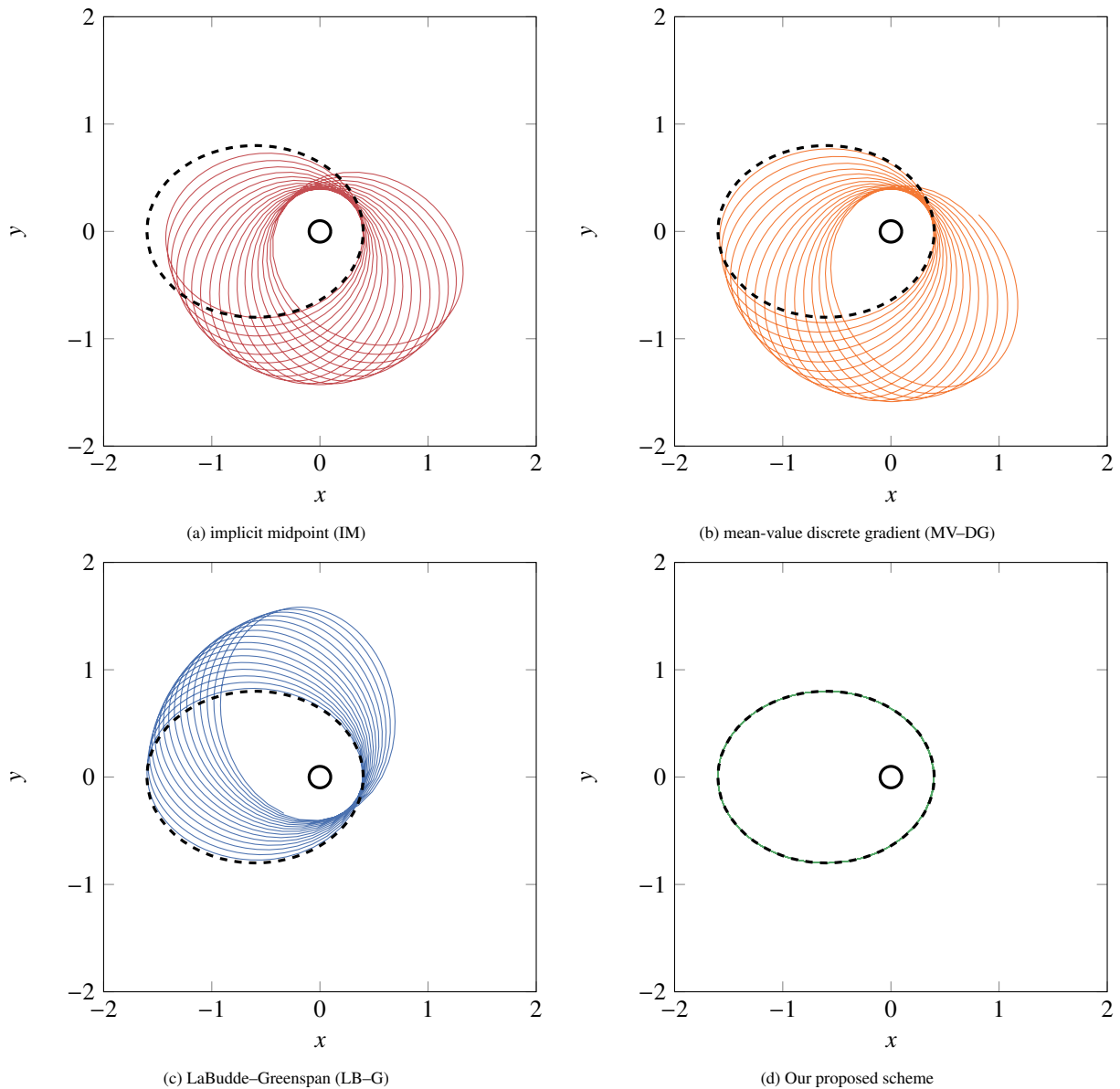

    \captionsetup[subfigure]{justification = centering}
    \centering

    \begin{subfigure}{0.5\textwidth}
        \centering



        \caption{Implicit midpoint (IM)}

        \label{fig:kepler_trajectories_implicit_midpoint}
    \end{subfigure}\begin{subfigure}{0.5\textwidth}
        \centering

        %

        \caption{AVF discrete gradient}

        \label{fig:kepler_trajectories_avf-dg}
    \end{subfigure}
    \begin{subfigure}{0.5\textwidth}
        \centering

        %

        \caption{LaBudde--Greenspan (LB--G)}

        \label{fig:kepler_trajectories_labudde_greenspan}
    \end{subfigure}\begin{subfigure}{0.5\textwidth}
        \centering

        %

        \caption{Our proposed scheme}

        \label{fig:kepler_trajectories_andrews_farrell}
    \end{subfigure}

    \caption{Trajectories of the Kepler problem. The exact solution is given by the dashed ellipse. All schemes are of the same order.}

    \label{fig:kepler_trajectories}
\end{figure}

To illustrate the qualitative benefits afforded by the fully conservative scheme \eqref{eq:kepler_avcpg}, we solve the Kepler problem with timestep $\Delta t_n = 0.1$ and final time $t = 100$, using various classical 1-stage, 2\textsuperscript{nd}-order implicit geometric integrators:
implicit midpoint (IM), the AVF discrete-gradient method of McLachlan, Quispel \& Robidoux \cite{McLachlan_Quispel_Robidoux_1999}, the LaBudde--Greenspan (LB--G) energy-- and angular momentum--conservative scheme \cite{LaBudde_Greenspan_1974}, and our scheme \eqref{eq:kepler_avcpg} at $s=1$.
\Cref{fig:kepler_trajectories} shows the resulting trajectories.

\begin{figure}[!ht]
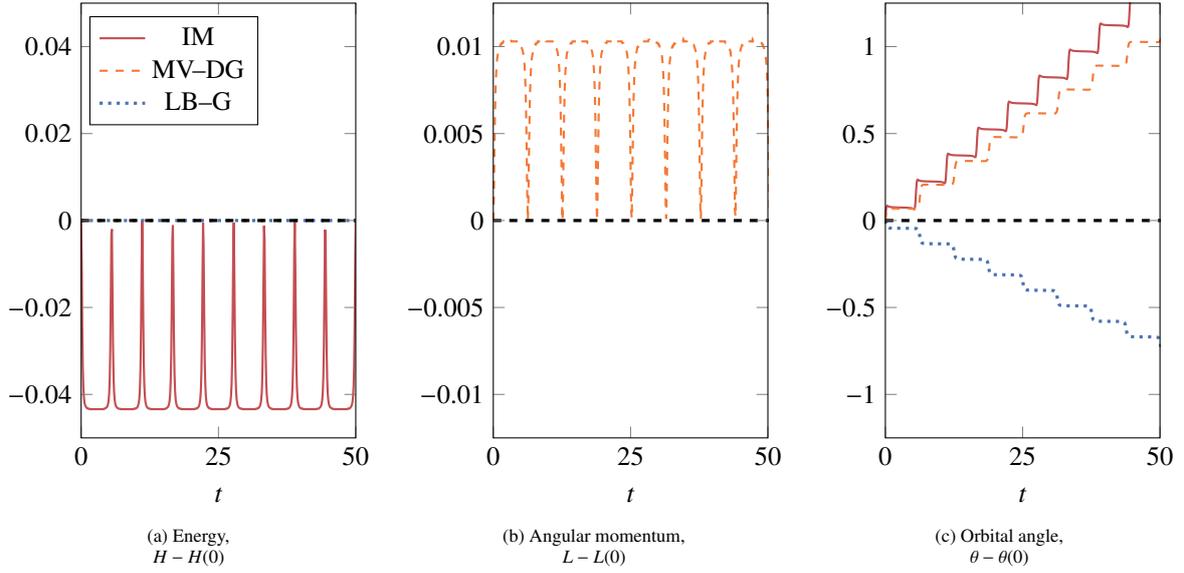

    \captionsetup[subfigure]{justification = centering}
    \centering

    \begin{subfigure}{0.319\textwidth}
        \centering



        \caption{Energy, \\ $H - H(0)$}
        \label{fig:kepler_invariants_im}
    \end{subfigure}
    \begin{subfigure}{0.319\textwidth}
        \centering

        %

        \caption{Angular momentum, \\ $L - L(0)$}
    \end{subfigure}
    \begin{subfigure}{0.319\textwidth}
        \centering

        %

        \caption{Orbital angle, \\ $\theta - \theta(0)$}
    \end{subfigure}

    \caption{Error in scalar invariants of the Kepler problem: $H$, $L$ and $\theta$.}

    \label{fig:kepler_invariants}
\end{figure}

In those cases where they are not conserved, \Cref{fig:kepler_invariants} shows the evolution of the invariants $H$, $L$, $\theta$ up to time $t = 50$, where $\theta \coloneqq \arg\bfA$.

As a symplectic method, implicit midpoint conserves the quadratic invariant $L$ (up to quadrature error, solver tolerances and machine precision) but neither $H$ nor $\theta$;
it therefore conserves neither the orbit shape nor its orientation, since trajectories in the Kepler problem should be precession-free.
The AVF discrete gradient scheme conserves $H$, but neither $L$ nor $\theta$.
The scheme of LaBudde \& Greenspan conserves $H$ and $L$ by design, and so conserves the orbit shape, but not its orientation $\theta$.
In contrast, our scheme \eqref{eq:kepler_avcpg} conserves all three invariants, thereby restricting the discrete solution to the same ellipse traced out by the exact solution.

These results illustrate the potential importance of conserving invariants in Poisson (and non-Poisson) systems:
while symplectic methods are likely preferable for e.g.~capturing the statistical behaviour of chaotic systems, conservative discretisations may give more physically reasonable results for individual trajectories at coarser timesteps.

\subsubsection{Convergence test}\label{sec:kepler_convergence}

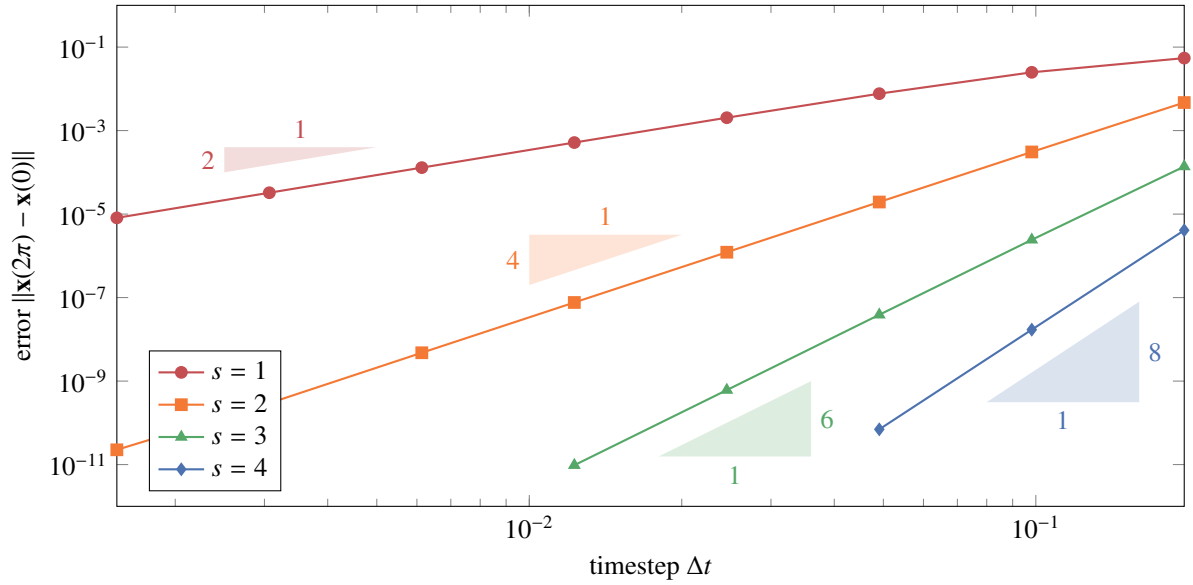
\begin{figure}[!ht]
    \begin{tikzpicture}
    \begin{axis}[
        xmode = log, xmin = 0.0015339807878856412, xmax = 0.19634954084936207, xlabel = {timestep $\Delta t$},
        ymode = log, ymin = 1e-12,                 ymax = 1e0,                 ylabel = {error $\|\bfx(2\pi) - \bfx(0)\|$},
        width = 0.95\textwidth, height = 0.5\textwidth,
        axis on top,
        legend pos=south west
    ]
        \addplot[seabornred, thick, mark=*] table {
    0.19634954084936207 0.05446609250338374
    0.09817477042468103 0.02488034338813408
    0.04908738521234052 0.007663864794629973
    0.02454369260617026 0.0020355302394051366
    0.01227184630308513 0.000517009496684054
    0.006135923151542565 0.00012977166781923005
    0.0030679615757712823 3.247555535211711e-05
    0.0015339807878856412 8.12092798909956e-06
};
         \addlegendentry{$s = 1$}
        \addplot[seabornorange, thick, mark=square*] table {
    0.19634954084936207 0.004698679690626697
    0.09817477042468103 0.0003077703925053552
    0.04908738521234052 1.9475125614559992e-05
    0.02454369260617026 1.2209853326754529e-06
    0.01227184630308513 7.637050019603455e-08
    0.006135923151542565 4.773867398444698e-09
    0.0030679615757712823 2.957160032767575e-10
    0.0015339807878856412 2.2602542964451686e-11
};
         \addlegendentry{$s = 2$}
        \addplot[seaborngreen, thick, mark=triangle*] table {
    0.19634954084936207 0.00013788954598268688
    0.09817477042468103 2.423282510617623e-06
    0.04908738521234052 3.8782267254967365e-08
    0.02454369260617026 6.096359394585306e-10
    0.01227184630308513 9.667441547620646e-12
};
         \addlegendentry{$s = 3$}
        \addplot[seabornblue, thick, mark=diamond*] table {
    0.19634954084936207 4.0969990019100395e-06
    0.09817477042468103 1.7065443181299557e-08
    0.04908738521234052 6.987954789059582e-11
};
         \addlegendentry{$s = 4$}

        \addplot[draw=none, fill=seabornred!20] table {
            0.005  0.0004
            0.0025 0.0004
            0.0025 0.0001
            0.005  0.0004
        };
        \addplot[draw=none] table {
            0.005  0.0004
            0.0025 0.0004
        } node[text=seabornred, pos=0.5, above] {$1$};
        \addplot[draw=none] table {
            0.0025 0.0004
            0.0025 0.0001
        } node[text=seabornred, pos=0.5, left] {$2$};
        \addplot[draw=none] table {
            0.0025 0.0001
            0.005  0.0004
        };
        
        \addplot[draw=none, fill=seabornorange!20] table {
            0.02 0.0000032
            0.01 0.0000032
            0.01 0.0000002
            0.02 0.0000032
        };
        \addplot[draw=none] table {
            0.02 0.0000032
            0.01 0.0000032
        } node[text=seabornorange, pos=0.5, above] {$1$};
        \addplot[draw=none] table {
            0.01 0.0000032
            0.01 0.0000002
        } node[text=seabornorange, pos=0.5, left] {$4$};
        \addplot[draw=none] table {
            0.01 0.0000002
            0.02 0.0000032
        };
        
        \addplot[draw=none, fill=seaborngreen!20] table {
            0.036 0.000000001
            0.036 0.000000000015625
            0.018 0.000000000015625
            0.036 0.000000001
        };
        \addplot[draw=none] table {
            0.036 0.000000001
            0.036 0.000000000015625
        } node[text=seaborngreen, pos=0.5, right] {$6$};
        \addplot[draw=none] table {
            0.036 0.000000000015625
            0.018 0.000000000015625
        } node[text=seaborngreen, pos=0.5, below] {$1$};
        \addplot[draw=none] table {
            0.018 0.000000000015625
            0.036 0.000000001
        };
        
        \addplot[draw=none, fill=seabornblue!20] table {
            0.16 0.00000008
            0.16 0.0000000003125
            0.08 0.0000000003125
            0.16 0.00000008
        };
        \addplot[draw=none] table {
            0.16 0.00000008
            0.16 0.0000000003125
        } node[text=seabornblue, pos=0.5, right] {$8$};
        \addplot[draw=none] table {
            0.16 0.0000000003125
            0.08 0.0000000003125
        } node[text=seabornblue, pos=0.5, below] {$1$};
        \addplot[draw=none] table {
            0.08 0.0000000003125
            0.16 0.00000008
        };
    \end{axis}
    \end{tikzpicture}

    \caption{Error in the position of the orbital body at $t = 2\pi$ for varying timesteps $\Delta t \in 2\pi \cdot 2^k$, $k \in \{-5, \dots, -12\}$ and stages $s \in \{1, \dots, 4\}$. The convergence curve for $s \in \{3, 4\}$ flattens out at smaller timesteps due to round-off error and solver tolerances. Triangles demonstrate observed convergence rates of $2s$.}

    \label{fig:kepler_convergence}
\end{figure}

\Cref{fig:kepler_convergence} shows the convergence of \eqref{eq:kepler_avcpg}, as measured by the error in the position of the orbital body after the true orbital period ($2\pi$ for these ICs) at varying timesteps $\Delta t_n$ and stages $s$.
As with the underlying Gau\ss{} scheme, we observe convergence with rate $2s$ before round-off errors dominate.

\subsection{Kovalevskaya problem}\label{sec:kovalevskaya}

As a further example, we consider the nondimensionalised Kovalevskaya top~\cite{Kovalevskaya_1889},
\begin{align}
      \dot{\bfn}  =  \bfn \times J\bfl,  \quad\quad
      \dot{\bfl}  =  \bfn \times \bfe_1 + \bfl \times J\bfl,
\end{align}
for $\bfn, \bfl : \bbR_+ \to \bbR^3$ representing the orientation vector (i.e.~the $z$-components of the principal axes) and the angular momentum (i.e.~the components of the angular momentum along those principal axes) respectively, $\bfe_1$ denoting the basis vector $(1, 0, 0)$, and $J$ denoting the matrix
\begin{equation}
    J  \coloneqq  \begin{pmatrix}
        1 & 0 & 0 \\
        0 & 1 & 0 \\
        0 & 0 & 2 \\
    \end{pmatrix}.
\end{equation}
Trajectories of this system have 4 invariants:
the energy $H \coloneqq \frac{1}{2}\bfl^\top J\bfl + \bfn\cdot\bfe_1$, the (square) norm of the orientation vector $\|\bfn\|^2$, the angular momentum in the $z$ direction $L = \bfl\cdot\bfn$, and the Kovalevskaya invariant $K = |\xi|^2$ where $\xi = (l_1 + il_2)^2 - 2(n_1 + in_2)$ ($i$ is the imaginary unit).
Note that while $H$, $\|\bfn\|^2$ and $\bfl\cdot\bfn$ are quadratic, $K$ is quartic.

Unlike the Kepler problem discussed above, it is not immediately clear from inspection how one might define an $\tilde{F} : \bbR^6 \to \Alt^5\bbR^6$ satisfying the conditions of \Cref{th:alternating_forms};
we therefore find such an $\tilde{F}$ using a construction similar to that used in the proof of \Cref{th:alternating_forms}.
Define the multilinear map $\tilde{G}((\bfn, \bfl)) : (\bbR^6)^5 \to \bbR$,
\begin{subequations}
\begin{equation}
    \tilde{G}\left(\begin{pmatrix} \bfn \\ \bfl \end{pmatrix}\right)\left[\begin{pmatrix} \bfa_1 \\ \bfb_1 \end{pmatrix}, \begin{pmatrix} \bfa_2 \\ \bfb_2 \end{pmatrix}, \begin{pmatrix} \bfa_3 \\ \bfb_3 \end{pmatrix}, \begin{pmatrix} \bfa_4 \\ \bfb_4 \end{pmatrix}, \begin{pmatrix} \bfm \\ \bfk \end{pmatrix}\right]
    \coloneqq
    \det[\bfb_1 \; \bfb_2 \; \bfb_3](\bfn\cdot\bfa_4)\left[
        \bfm^\top(\bfn \times J\bfl) + \bfk^\top(\bfn \times \bfe_1) + \bfk^\top(\bfl \times J\bfl)
    \right].
\end{equation}
Considering the alternatisation $\Alt\tilde{G}((\bfn, \bfl)) \in \Alt^5\bbR^6$, if we apply $\Alt\tilde{G}((\bfn, \bfl))$ to the gradients of the invariants $H$, $K$, $L$, $\frac{1}{2}\|\bfn\|^2$, we see
\begin{equation}
    \Alt\tilde{G}\left(\begin{pmatrix} \bfn \\ \bfl \end{pmatrix}\right)\left[\begin{pmatrix} \bfe_1 \\ J\bfl \end{pmatrix}, \begin{pmatrix} \nabla_\bfn K \\ \nabla_\bfl K \end{pmatrix}, \begin{pmatrix} \bfl \\ \bfn \end{pmatrix}, \begin{pmatrix} \bfn \\ \bfzero \end{pmatrix}, \begin{pmatrix} \bfm \\ \bfk \end{pmatrix}\right]
    =
    6\det[J\bfl \; \nabla_\bfl K \; \bfn]\|\bfn\|^2\left[
        \bfm^\top(\bfn \times J\bfl) + \bfk^\top(\bfn \times \bfe_1) + \bfk^\top(\bfl \times J\bfl)
    \right].
\end{equation}
We therefore define $\tilde{F}((\bfn, \bfl)) \in \Alt^5\bbR^6$,
\begin{equation}
\label{eq:F_kovalevskaya}
    \tilde{F}\left(\begin{pmatrix} \bfn \\ \bfl \end{pmatrix}\right)\left[\begin{pmatrix} \bfa_1 \\ \bfb_1 \end{pmatrix}, \begin{pmatrix} \bfa_2 \\ \bfb_2 \end{pmatrix}, \begin{pmatrix} \bfa_3 \\ \bfb_3 \end{pmatrix}, \begin{pmatrix} \bfa_4 \\ \bfb_4 \end{pmatrix}, \begin{pmatrix} \bfm \\ \bfk \end{pmatrix}\right]
    \coloneqq
    \frac{1}{6\det[J\bfl \; \nabla_\bfl K \; \bfn]\|\bfn\|^2}\Alt\tilde{G}\left(\begin{pmatrix} \bfn \\ \bfl \end{pmatrix}\right)\left[\begin{pmatrix} \bfa_1 \\ \bfb_1 \end{pmatrix}, \begin{pmatrix} \bfa_2 \\ \bfb_2 \end{pmatrix}, \begin{pmatrix} \bfa_3 \\ \bfb_3 \end{pmatrix}, \begin{pmatrix} \bfa_4 \\ \bfb_4 \end{pmatrix}, \begin{pmatrix} \bfm \\ \bfk \end{pmatrix}\right].
\end{equation}
\end{subequations}
This satisfies \eqref{eq:alternating_form_coincidence}.
We may then use this $\tilde{F}$ in \eqref{eq:integrable_avcpg} to define a fully conservative integrator for the Kovalevskaya top.

\subsubsection{Numerical test}\label{sec:kovalevskaya_test}

\revised{%
To test our proposed conservative integrator, we consider three discretisations of the Kovalevskaya top:
(i)~implicit midpoint, (ii)~implicit midpoint with projection, i.e.~projecting back onto the invariant manifold after each timestep \cite[Sec.~IV.4]{Hairer_Lubich_Wanner_2006}, and (iii)~our proposed} fully conservative modification of 1-stage continuous Petrov--Galerkin \eqref{eq:integrable_avcpg} (i.e.~setting $\calI_n$ to be the exact integral) using $\tilde{F}$ as in \eqref{eq:F_kovalevskaya}.
We use the same ICs $\bfn(0) = (0.8, 0.6, 0)$, $\bfl(0) = (2, 0, 0.2)$ in each case, and a uniform timestep $\Delta t_n = 0.1$ up to a final time $t = 300$.

\begin{figure}[!ht]
    \captionsetup[subfigure]{justification = centering}
    \centering

    \begin{subfigure}{0.5\textwidth}
        \raggedright
        \includegraphics[width = 0.7\textwidth]{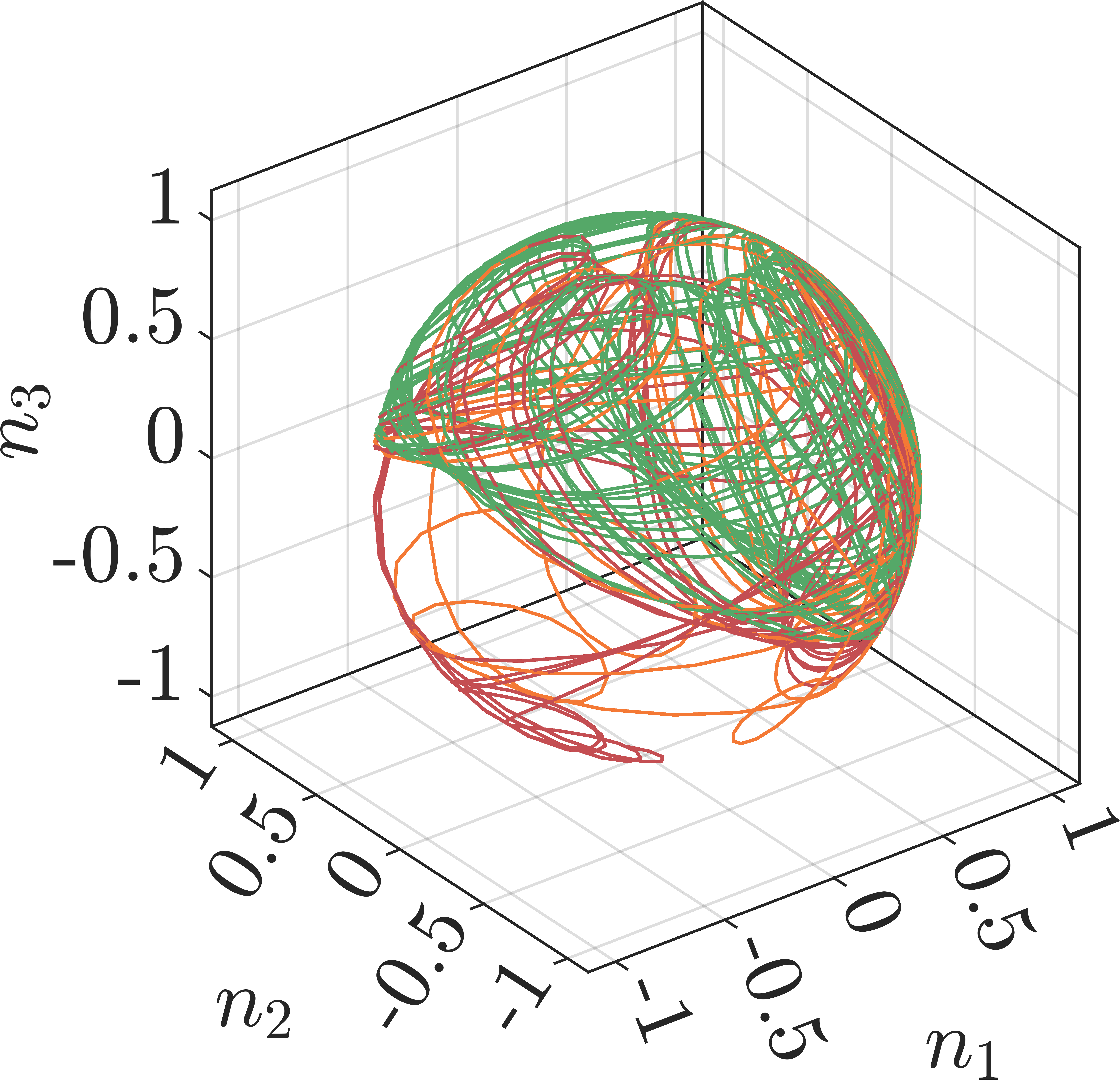}

        \raggedleft
        \includegraphics[width = 0.7\textwidth]{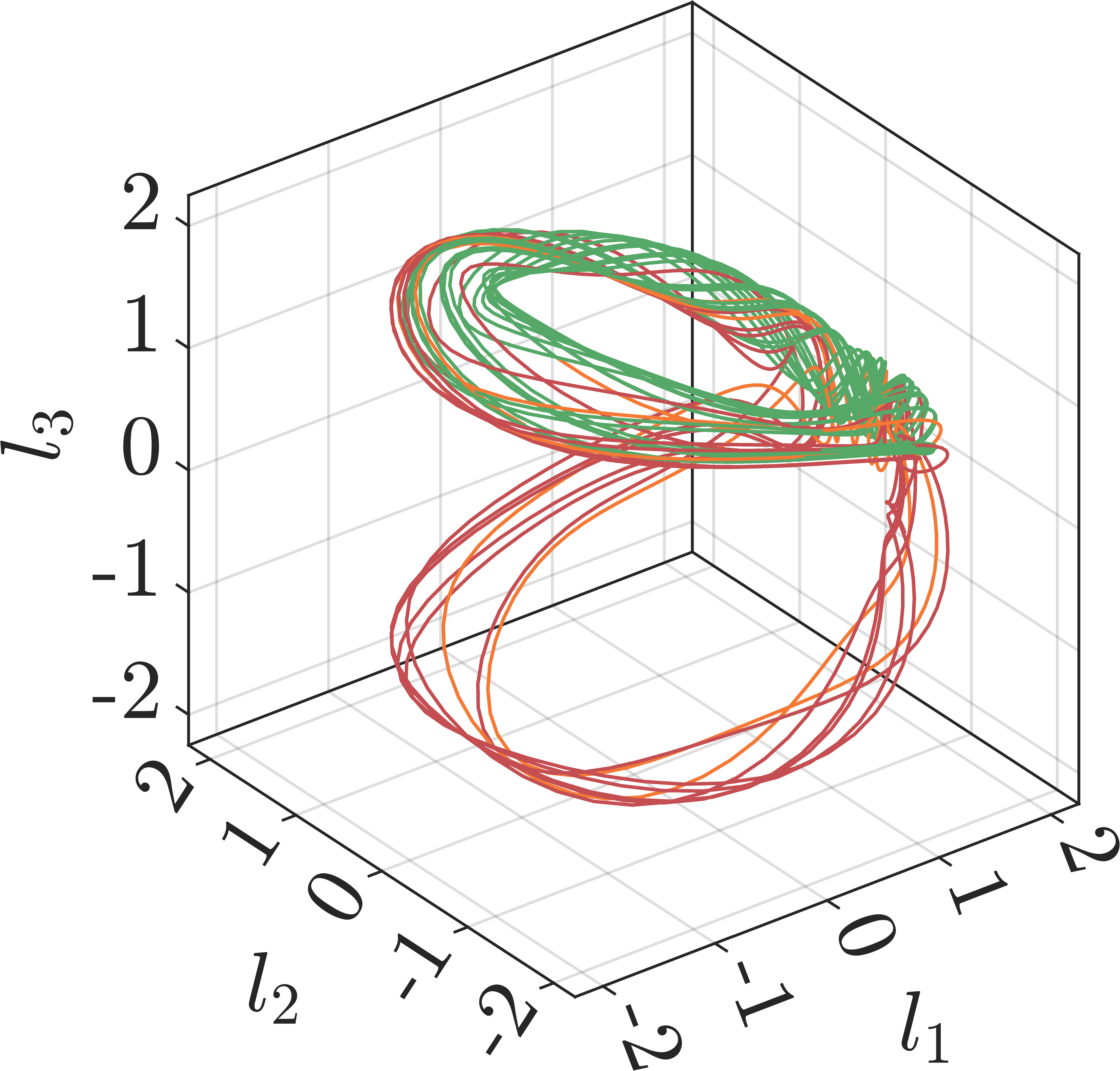}

        \centering
        \caption{Implicit midpoint}

        \label{fig:kovalevskaya_trajectories_implicit_midpoint}
    \end{subfigure}\begin{subfigure}{0.5\textwidth}
        \raggedright
        \includegraphics[width = 0.7\textwidth]{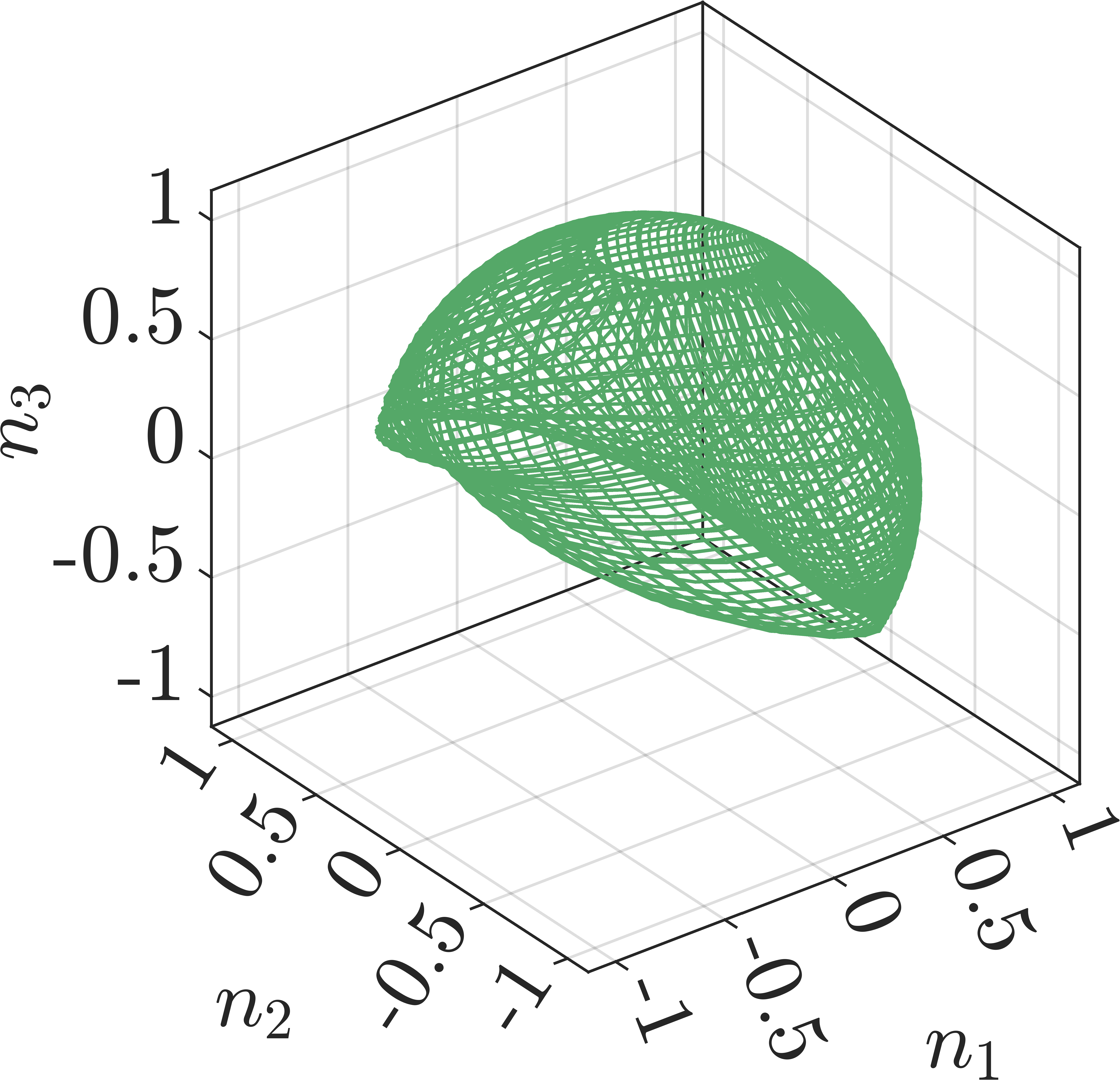}

        \raggedleft
        \includegraphics[width = 0.7\textwidth]{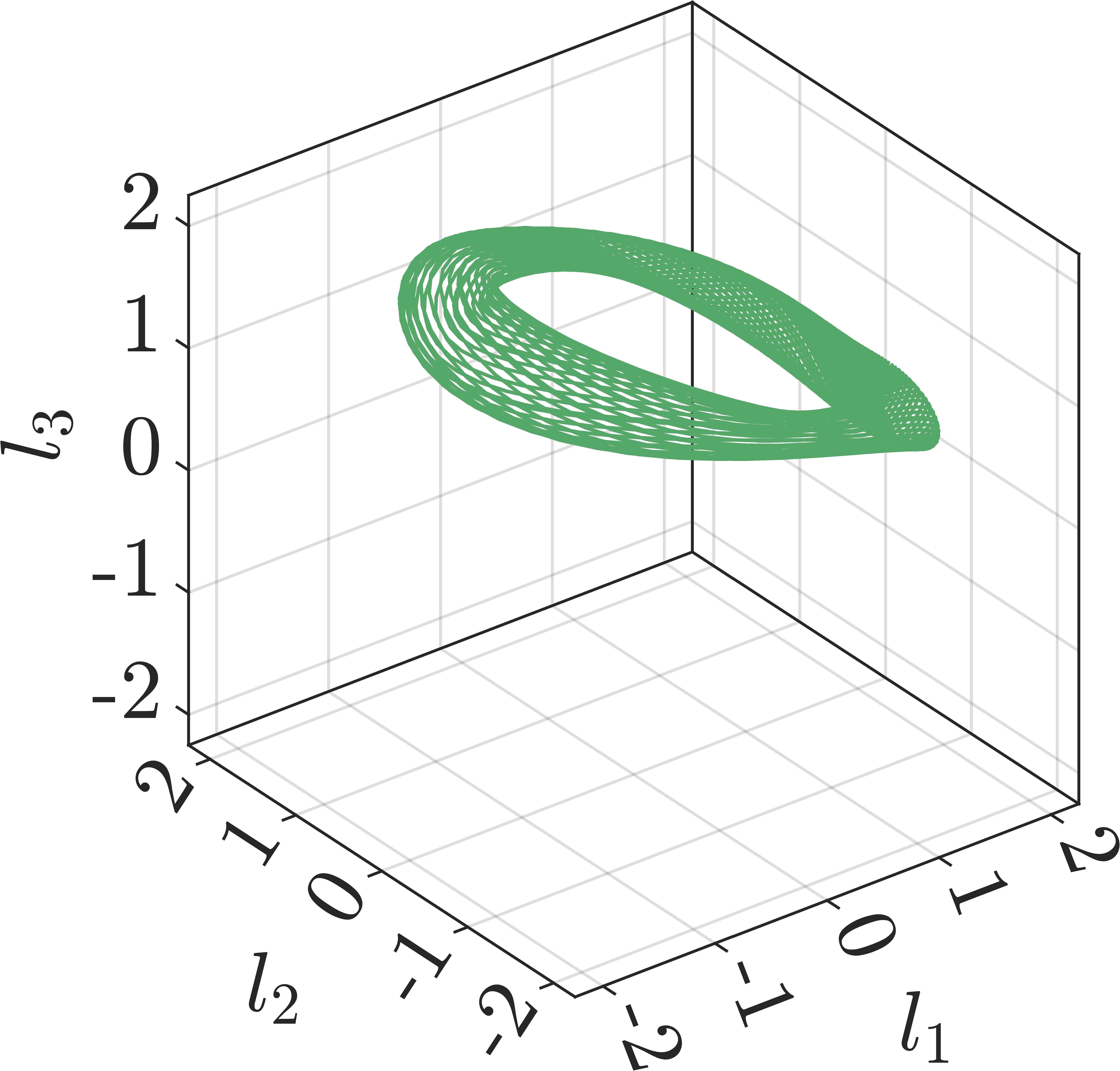}

        \centering
        \caption{Our scheme}

        \label{fig:kovalevskaya_trajectories_andrews_farrell}
    \end{subfigure}

    \caption{Trajectories in $\bfn, \bfl$ of the Kovalevskaya top, with implicit midpoint (left) and our proposed scheme (right).}

    \label{fig:kovalevskaya_trajectories}
\end{figure}

\begin{figure}[!ht]
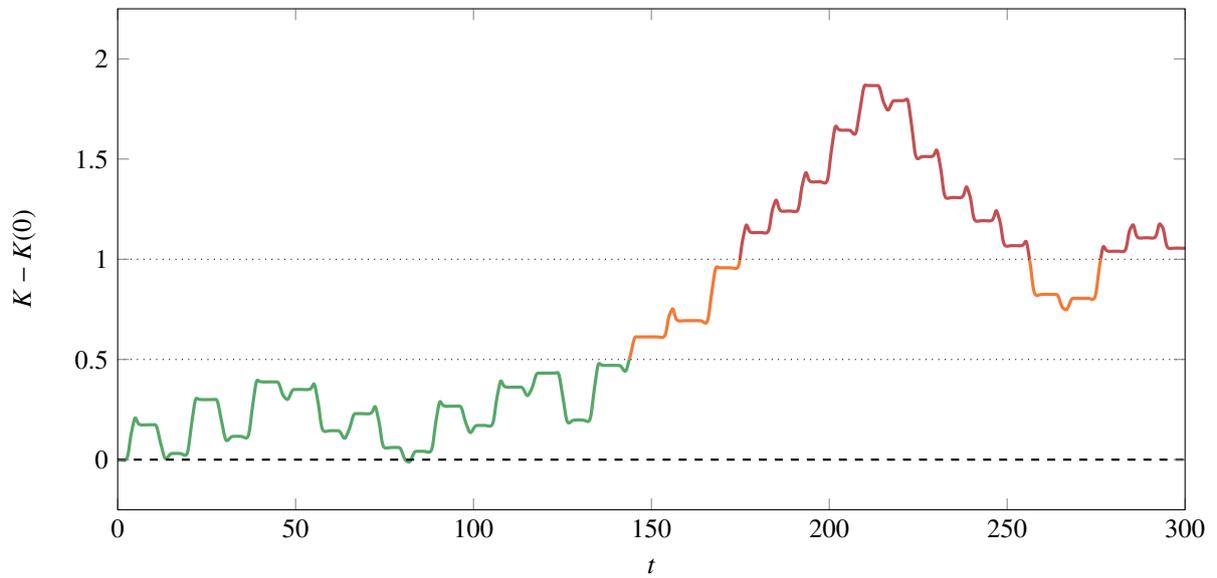

    \centering



    \caption{Error $K - K(0)$ within the implicit midpoint simulation of the Kovalevskaya top.}

    \label{fig:kovalevskaya_invariant}
\end{figure}

\revised{%
\Cref{fig:kovalevskaya_trajectories} shows discrete trajectories for the implicit midpoint discretisation alongside our proposed scheme.
Despite the projection method sacrificing time symmetry (something which is retained by our scheme), the numerical solution differs in the $\ell^2$ norm from ours only by at most $2.7\times10^{-3}$ in $\bfn$ and $1.9\times10^{-3}$ in $\bfl$ over the duration.
We therefore do not plot the associated trajectory as it is visually indistinguishable.
}
\Cref{fig:kovalevskaya_invariant} shows the evolution and drift of the Kovalevskaya invariant $K$ within the implicit midpoint scheme.
In both figures, colouring indicates error in the Kovalevskaya invariant $K$:
green for $|K - K(0)| \le \frac{1}{2}$, orange for $|K - K(0)| \in (\frac{1}{2}, 1]$, red for $|K - K(0)| > 1$.

All invariants, including $K$, are conserved by the trajectory of our scheme (up to quadrature error, solver tolerances and machine precision).
\revised{%
The quadratic invariants $H$, $\|\bfn\|^2$ and $L$ are all conserved by the symplectic implicit midpoint scheme \cite[Sec.~IV.2]{Hairer_Lubich_Wanner_2006}.
Being quartic, however, $K$ is not conserved by implicit midpoint, and is the only invariant it fails to preserve;
}
we see that the resulting drift in $K$ allows the numerical simulation to admit unphysical trajectories after a sufficient duration (approximately 18 rotations of the top for these ICs).
\revised{%
This demonstrates directly the value of conserving $K$, and more generally that conserving only those invariants that happen to be quadratic is not always sufficient.
}

\section{GENERIC ODEs}\label{sec:generic_odes}

We now turn to systems of ODEs in the GENERIC formalism \cite{Grmela_Ottinger_1997, Ottinger_Grmela_1997}, which couple reversible (Poisson) and irreversible (gradient-descent) dynamics in a thermodynamically compatible way.
The general GENERIC ODE has the form
\begin{equation}\label{eq:generic_ode}
    \dot{\bfx} = B(\bfx)\nabla E(\bfx) + D(\bfx)\nabla S(\bfx),
\end{equation}
where $B(\bfx), D(\bfx) \in \bbR^{d\times d}$ are respectively skew-symmetric and positive-semidefinite (the Poisson and friction matrices), and $E, S : \bbR^d \to \bbR$ denote the energy and entropy.
This augments the general Poisson system \eqref{eq:poisson} with an additional dissipative contribution.
Thermodynamic consistency imposes the orthogonality (degeneracy) conditions
\begin{equation}\label{eq:generic_ode_compatibility}
    \nabla S(\bfx)^{\top} B(\bfx)  =  \bfzero^\top, \qquad
    \nabla E(\bfx)^{\top} D(\bfx)  =  \bfzero^\top.
\end{equation}
We see these conditions ensure energy conservation ($\tfrac{\rmd}{\rmd t}E(\bfx) = 0$) and entropy generation ($\tfrac{\rmd}{\rmd t}S(\bfx) \ge 0$) along exact trajectories by testing \eqref{eq:generic_ode} against $\nabla E$ and $\nabla S$.

\begin{remark}\label{rem:metriplectic}
    \revised{%
    \emph{Metriplectic systems} \cite{Morrison_1986} represent a generalisation of Poisson systems \eqref{eq:poisson}, incorporating a second positive-semidefinite bracket alongside the Poisson bracket\footnote{
        \revised{The introduction of dissipative effects into Poisson systems in this bracket form dates to the 1984 works of Grmela \cite{Grmela_1984}, Kaufman \cite{Kaufman_1984} and Morrison \cite{Morrison_1984a, Morrison_1984b}, the last interpreting the dissipation as an entropy generation.
        The terminology of \emph{metriplectic} systems, and the reading of the entropy as a Casimir of the underlying Poisson system with the dissipative term induced by a metric, followed in \cite{Morrison_1986};
        it remains the more common usage in the plasma physics and Hamiltonian mechanics communities, whereas the title of GENERIC \cite{Grmela_Ottinger_1997, Ottinger_Grmela_1997}, attached to the tensorial formulation, came later.
        In the metriplectic formulation, $E$ would typically be denoted $H$; however, we use $E$ here for consistency.}
    }:
    \begin{equation}\label{eq:metriplectic}
        \dot{\bfx} = \{\bfx, E\} + (\bfx, S\!).
    \end{equation}
    The brackets must satisfy the compatibility conditions\footnote{
        \revised{For the former condition, $S$ is equivalently a Casimir of the original Poisson system.}
    }
    \begin{equation}\label{eq:metriplectic_compatibility}
        \{S, \cdot\}  =  0, \qquad
         (E, \cdot)   =  0,
    \end{equation}
    guaranteeing energy conservation and entropy generation.
    The matrix-formulated GENERIC system \eqref{eq:generic_ode} relates to the bracket-formulated metriplectic system \eqref{eq:metriplectic}, just as the matrix formulation \eqref{eq:poisson_tensor} of the Poisson system relates to the bracket formulation \eqref{eq:poisson_bracket}.
    As in \Cref{rem:nambu}, it is helpful for our purposes to consider the matrix formulation, since it is written in terms of the gradients $\nabla E$ and $\nabla S$, and so offers a natural place to introduce the auxiliary variables.
    }
\end{remark}

To preserve these thermodynamic properties on the discrete level, the strategy in the construction of \eqref{eq:odecpgmod} and \eqref{eq:integrable_avcpg} then suggests (i)~introducing auxiliary variables $\widetilde{\nabla E}$ and $\widetilde{\nabla S}$ approximating $\nabla E$ and $\nabla S$ respectively, and (ii)~introducing these in the right-hand side of \eqref{eq:generic_ode}.
It is, however, insufficient to substitute these auxiliary variables only for the gradients seen in \eqref{eq:generic_ode}, as this will fail to reproduce the orthogonality conditions \eqref{eq:generic_ode_compatibility}.
We therefore assume access to certain consistent extensions of $B$ and $D$ that can be evaluated on these auxiliary variables.

\begin{assumption}[Characterisation of GENERIC matrix compatibility]\label{ass:generic_matrices}
    Assume the existence of $\tilde{B}, \tilde{D} : (\bbR^d)^2 \to \bbR^{d\times d}$ such that the following hold:
    \begin{subequations}
    \begin{enumerate}
        \item Consistency with $B, D$:
        for all $\bfx \in \bbR^d$,
        \begin{align}\label{eq:generic_matrix_coincidence}
            \tilde{B}(\bfx, \nabla S(\bfx)) = B(\bfx), \qquad
            \tilde{D}(\bfx, \nabla E(\bfx)) = D(\bfx).
        \end{align}
        \item Skew-symmetry and positive-semidefiniteness:
        for all $\bfx, \widetilde{\nabla E}, \widetilde{\nabla S} \in \bbR^d$, $\tilde{B}(\bfx, \widetilde{\nabla S})$ is skew-symmetric, and $\tilde{D}(\bfx, \widetilde{\nabla E})$ is positive-semidefinite.
        \item Preservation of compatibility \eqref{eq:generic_ode_compatibility}:
        for all $\bfx, \widetilde{\nabla E}, \widetilde{\nabla S} \in \bbR^d$,
        \begin{align}\label{eq:generic_matrix_compatibility}
            \widetilde{\nabla S}^{\top}\tilde{B}(\bfx, \widetilde{\nabla S})  =  \bfzero^\top, \qquad
            \widetilde{\nabla E}^{\top}\tilde{D}(\bfx, \widetilde{\nabla E})  =  \bfzero^\top.
        \end{align}
    \end{enumerate}
    \end{subequations}
\end{assumption}

\begin{remark}\label{rem:ode_scope}
    While we cannot guarantee \Cref{ass:generic_matrices} always holds in full generality, it holds in the typical cases we consider below, and suffices for our discrete stability results.
    \revised{%
    The orthogonality conditions \eqref{eq:generic_ode_compatibility} typically hold as natural consequences of how the matrices $B$ and $D$ are constructed in terms of $E$ and $S$, or more specifically their gradients $\nabla E$ and $\nabla S$.
    One can then obtain $\tilde{B}$ and $\tilde{D}$ by identifying each occurrence of $\nabla E$ and $\nabla S$ in the definition of $B$ and $D$, and replacing them with the corresponding auxiliary variable.
    We perform this natural substitution to define $\tilde{D}$ (a function of $\widetilde{\nabla E}$) in the thermodynamic engine example below (\Cref{sec:engine}).
    A similar argument applies in the PDE case below (\Cref{sec:generic_pdes}), where $B$ and $D$ constitute skew-symmetric and positive-semidefinite operators, and a similar compatibility assumption is required (\Cref{ass:generic_operators}).
    Theoretically, if no natural $\tilde{B}$ and $\tilde{D}$ could be found satisfying \Cref{ass:generic_matrices}, one option could be to construct them directly, projecting $B$ and $D$ onto the orthogonal complements of $\widetilde{\nabla S}$ and $\widetilde{\nabla E}$ respectively.
    Care would have to be taken, however, to ensure that the skew-symmetry and positive-semidefiniteness are preserved, and we have not yet found this to be necessary in practice.
    }
\end{remark}

With this in hand, the GENERIC ODE \eqref{eq:generic_ode} can be discretised in an energy- and entropy-stable manner by defining auxiliary variables $\widetilde{\nabla E}, \widetilde{\nabla S} \in \dot{\bbX}_n$, introducing these into the right-hand side, and replacing the Poisson and friction matrices with their consistent extensions evaluated at these auxiliary variables:
find $(\bfx,(\widetilde{\nabla E},\widetilde{\nabla S})) \in \bbX_n \times (\dot{\bbX}_n)^2$ such that, for all $(\bfy,(\bfy_E,\bfy_S)) \in \dot{\bbX}_n \times (\dot{\bbX}_n)^2$,
\begin{subequations}\label{eq:generic_ode_avcpg}
\begin{align}
    \calI_n\left[\bfy^{\top}\dot{\bfx}\right]  &= \calI_n\left[\bfy^{\top}\tilde{B}(\bfx, \widetilde{\nabla S})\widetilde{\nabla E} + \bfy^{\top}\tilde{D}(\bfx, \widetilde{\nabla E})\widetilde{\nabla S}\right], \label{eq:generic_ode_avcpg_a} \\
    \calI_n\left[\widetilde{\nabla E}^{\top}\bfy_E\right]  &= \int_{T_n} \nabla E(\bfx)^{\top}\bfy_E, \\
    \calI_n\left[\widetilde{\nabla S}^{\top}\bfy_S\right]  &= \int_{T_n} \nabla S(\bfx)^{\top}\bfy_S.
\end{align}
\end{subequations}

\begin{theorem}[Energy and entropy stability of \eqref{eq:generic_ode_avcpg}]\label{th:generic_ode_sp}
    Assuming solutions exist, the integrator \eqref{eq:generic_ode_avcpg} is energy- and entropy-stable, with $E(\bfx(t_{n+1})) = E(\bfx(t_n))$ and $S(\bfx(t_{n+1})) \ge S(\bfx(t_n))$.
\end{theorem}

\begin{proof}
    By considering respectively $\bfy_E = \dot{\bfx}$, $\bfy_S = \dot{\bfx}$ and $\bfy = \widetilde{\nabla E}$, $\bfy = \widetilde{\nabla S}$ in \eqref{eq:generic_ode_avcpg},
    \begin{subequations}
    \begin{align}
        E(\bfx(t_{n+1})) &- E(\bfx(t_n))  &
            S(\bfx(t_{n+1})) &- S(\bfx(t_n))  \notag  \\
        &= \int_{T_n} \dot{E}  &
            &= \int_{T_n} \dot{S}  \\
        &= \int_{T_n} \nabla E(\bfx)^\top \dot{\bfx}  &
            &= \int_{T_n} \nabla S(\bfx)^\top \dot{\bfx}  \\
        &= \calI_n\left[\widetilde{\nabla E}^\top \dot{\bfx}\right]  &
            &= \calI_n\left[\widetilde{\nabla S}^\top \dot{\bfx}\right]  \\
        &= \calI_n\left[\begin{aligned}
            &\widetilde{\nabla E}^{\top}\tilde{B}(\bfx, \widetilde{\nabla S})\widetilde{\nabla E}  \\
            &\qquad\quad+ \widetilde{\nabla E}^{\top}\tilde{D}(\bfx, \widetilde{\nabla E})\widetilde{\nabla S}
        \end{aligned}\right]  &
            &= \calI_n\left[\begin{aligned}
                &\widetilde{\nabla S}^{\top}\tilde{B}(\bfx, \widetilde{\nabla S})\widetilde{\nabla E}  \\
                &\qquad\qquad+ \widetilde{\nabla S}^{\top}\tilde{D}(\bfx, \widetilde{\nabla E})\widetilde{\nabla S}
            \end{aligned}\right]  \\
        &= 0,  &
            &\ge  0,
    \end{align}
    \end{subequations}
    where the final equality and inequality hold by \Cref{ass:generic_matrices}.
\end{proof}

\begin{remark}\label{rem:generic_ode_cost}
    \revised{%
    As with \eqref{eq:odecpgmod}, the auxiliary variables introduced in \eqref{eq:generic_ode_avcpg} can be eliminated, such that no additional unknowns enter the nonlinear system;
    the cost per timestep is that of (i)~the underlying Runge--Kutta method, together with (ii)~whatever additional quadrature the exact integrals demand.
    This is in contrast to the PDE setting of \Cref{sec:generic_pdes} below, where such an elimination is not in general available (see \Cref{rem:generic_pde_cost}).
    Despite this, our implementation of the engine discretisation \eqref{eq:engine_avcpg} below retains the proposed auxiliary temperature variables $(\tilde{T}_c)_{c=1}^C$ (components of the discrete gradient $\widetilde{\nabla E}$) as unknowns, as we find the increased computational cost to be minimal\footnote{
        \revised{%
        The algebraic systems are small enough that the linear solves account for under $0.5\%$ of the time spent, essentially all of which goes on assembling the residual and Jacobian.
        This cost is essentially unchanged by the introduction of $(\tilde{T}_c)_{c=1}^C$:
        across stages $s \in \{1, 2, 3\}$ we measure between $22$ and $24$ milliseconds of CPU time per timestep under either scheme, with neither scheme consistently faster than the other.
        Nor does our proposed discretisation require more Newton iterations per timestep than the Gau\ss{} method, taking on average $4.2$ against $4.3$.
        }
    }.}
\end{remark}

\begin{remark}\label{rem:stiff_problems}
    \revised{%
    The conservation and dissipation properties of \Cref{th:generic_ode_sp} hold whenever a discrete solution exists, regardless of the timestep.
    However, the accuracy and behaviour of \eqref{eq:generic_ode_avcpg} are inherited from the underlying continuous Petrov--Galerkin time discretisation.
    Like Gau\ss{} methods, such discretisations are A-stable but not L-stable;
    stiff modes that should decay quickly are expected to persist as spurious oscillations.
    This is a concern for GENERIC systems in which the friction matrix $D$ drives entropy production on a timescale short by comparison with the timestep $\Delta t_n$.
    The schemes proposed here would therefore typically be inadvisable in such stiff settings, where an energy-conserving L-stable discretisation might be preferred.
    As ever, the appropriate choice of method depends on the problem setting.
    }
\end{remark}

\subsection{Internal combustion engine}\label{sec:engine}

Inspired by the classical thermodynamic systems considered by e.g.~\"Ottinger \cite[Ex.~3]{Ottinger_2005} or Gay-Balmaz \& Yoshimura \cite[Sec.~3.1]{Gay-Balmaz_Yoshimura_2017}, we illustrate the GENERIC ODE discretisation \eqref{eq:generic_ode_avcpg} with a simple nondimensionalised model of a $C$-cylinder internal combustion engine, exchanging heat with an isothermal environment.
The engine in our model is unfired, i.e.~no combustion events are occurring, and coasting, i.e.~spinning on its own inertia alone;
together with the environmental heat exchange, this ensures the system is dissipative, fitting into the GENERIC formalism.

Let $\theta$ denote the engine phase, measured through an angular displacement, and $\omega$ its rate of change.
For cylinders $c=1,\dots,C$ we denote by $P_c$, $T_c$, $U_c$ and $S_c$ the pressure, temperature, energy and entropy of the working fluid in cylinder $c$;
these satisfy constitutive relations associated with the medium, completed by the volume $V_c$ given as a function of $\theta$ as
\begin{equation}\label{eq:piston_volume}
    V_c \coloneqq V_p - \cos\left(\theta - \frac{2\pi c}{C}\right),
\end{equation}
where $V_p > 1$ is the average volume across the pistons.
We denote the temperature, energy and entropy of the surrounding environment by $T_0$, $U_0$ and $S_0$, again coupled by constitutive relations, with $T_0$ assumed constant.
The engine model is then given by
\begin{subequations}\label{eq:engine_system}
\begin{align}
    \dot{\theta} &= \omega, &
    \dot{S}_c &= \frac{T_0 - T_c}{T_c}, \\
    \dot{\omega} &= \sum_{c=1}^C P_c \sin\left(\theta - \frac{2\pi c}{C}\right), &
    \dot{S}_0 &= \sum_{c=1}^C\frac{T_c - T_0}{T_0}.
\end{align}
\end{subequations}
We collect the state as $\bfx \coloneqq (\theta,\omega,(S_c)_{c=1}^C,S_0)$, and define the total energy and entropy
\begin{equation}
    E(\bfx) \coloneqq \frac{1}{2}\omega^2 + \sum_{c=1}^C U_c + U_0, \qquad
    S(\bfx) \coloneqq \sum_{c=1}^C S_c + S_0;
\end{equation}
as above, the energies $(U_c)$ and $U_0$ may be written as functions of the state $\bfx$ through the constitutive relations of the medium, with the piston volumes $(V_c)$ given by \eqref{eq:piston_volume}. 
From the fundamental thermodynamic relations $\rmd U_c = T_c\,\rmd S_c - P_c\,\rmd V_c$ and $\rmd U_0 = T_0\,\rmd S_0$, we obtain
\begin{equation}
    \nabla E(\bfx) = \begin{pmatrix}
        -\sum_{c=1}^C P_c\sin\left(\theta - \frac{2\pi c}{C}\right) \\
        \omega \\
        (T_c)_{c=1}^C \\
        T_0
    \end{pmatrix},
    \qquad
    \nabla S(\bfx) = \begin{pmatrix} 0 \\ 0 \\ \bfone \\ 1 \end{pmatrix}.
\end{equation}
The GENERIC representation \eqref{eq:generic_ode} of the combustion engine model \eqref{eq:engine_system} follows with a Poisson matrix $B$ equal to the canonical $2\times2$ skew block on $(\theta,\omega)$ and zeros elsewhere, and a positive-semidefinite friction matrix $D$ encoding thermal relaxation between each cylinder and the environment:
\begin{equation}
    B(\bfx) = \begin{pmatrix}
        0 & 1 & 0 & \cdots & 0 & 0 \\
        -1 & 0 & 0 & \cdots & 0 & 0 \\
        0 & 0 & 0 & \cdots & 0 & 0 \\
        \vdots & \vdots & \vdots & \ddots & \vdots & \vdots \\
        0 & 0 & 0 & \cdots & 0 & 0 \\
        0 & 0 & 0 & \cdots & 0 & 0
    \end{pmatrix},
    \qquad
    D(\bfx) = \begin{pmatrix}
        0 & 0 & 0 & \cdots & 0 & 0 \\
        0 & 0 & 0 & \cdots & 0 & 0 \\
        0 & 0 & \tfrac{T_0}{T_1} & & 0 & -1 \\
        \vdots & \vdots & & \ddots & & \vdots \\
        0 & 0 & 0 & & \tfrac{T_0}{T_C} & -1 \\
        0 & 0 & -1 & \cdots & -1 & \sum\limits_{c=1}^C \tfrac{T_c}{T_0}
    \end{pmatrix}.
\end{equation}
The degeneracy conditions \eqref{eq:generic_ode_compatibility} hold:
$\nabla S(\bfx)^{\top}B(\bfx) = \bfzero^\top$ by sparsity, and $\nabla E(\bfx)^{\top}D(\bfx) = \bfzero^\top$ by the structure of the thermal block.

For the discrete scheme \eqref{eq:generic_ode_avcpg} we require $\tilde{B}$ and $\tilde{D}$ satisfying \Cref{ass:generic_matrices}.
For the former, we may simply take $\tilde{B} = B$ as $B$ is independent of the state $\bfx$;
for the latter, a consistent $\tilde{D}$ may be defined by replacing each temperature $T_c$ with the corresponding components of $\widetilde{\nabla E}$:
\begin{equation}
    \tilde{D}(\bfx, \widetilde{\nabla E}) \coloneqq \begin{pmatrix}
        0 & 0 & 0 & \cdots & 0 & 0 \\
        0 & 0 & 0 & \cdots & 0 & 0 \\
        0 & 0 & \tfrac{T_0}{\widetilde{\partial_{S_1}E}} & & 0 & -1 \\
        \vdots & \vdots & & \ddots & & \vdots \\
        0 & 0 & 0 & & \tfrac{T_0}{\widetilde{\partial_{S_C}E}} & -1 \\
        0 & 0 & -1 & \cdots & -1 & \sum\limits_{c=1}^C \tfrac{\widetilde{\partial_{S_c}E}}{T_0}
    \end{pmatrix}.
\end{equation}

Neglecting auxiliary variables that turn out to be constant, the energy- and entropy-stable discretisation given by \eqref{eq:generic_ode_avcpg} is then as follows:
find $((\theta,\omega,(S_c),S_0), (\tilde{P},\tilde{\omega},(\tilde{T}_c)_{c=1}^C)) \in \bbX_n^{C+3} \times \dot{\bbX}_n^{C+2}$ such that, for all $((\eta,\psi,(R_c)_{c=1}^C,R_0),(\tilde{\psi},\tilde{\eta},(\tilde{W}_c)_{c=1}^C)) \in \dot{\bbX}_n^{C+3} \times \dot{\bbX}_n^{C+2}$,
\begin{subequations}\label{eq:engine_avcpg}
\begin{align}
    \calI_n[\dot{\theta}\eta] &= \calI_n[\tilde{\omega}\eta], \label{eq:engine_avcpg_a} \\
    \calI_n[\dot{\omega}\psi] &= - \calI_n[\tilde{P}\psi], \label{eq:engine_avcpg_b} \\
    \calI_n[\dot{S}_cR_c] &= \calI_n\left[\frac{T_0-\tilde{T}_c}{\tilde{T}_c}R_c\right],  \qquad  c = 1, \dots, C, \\
    \calI_n[\dot{S}_0R_0] &= \calI_n\left[\sum_{c=1}^C \frac{\tilde{T}_c - T_0}{T_0}R_0\right], \\
    \calI_n[\tilde{P}\tilde{\psi}] &= \int_{T_n} \sum_{c=1}^C P_c\sin\left(\theta - \frac{2\pi c}{C}\right)\tilde{\psi}, \label{eq:engine_avcpg_e} \\
    \calI_n[\tilde{\omega}\tilde{\eta}] &= \int_{T_n}\omega\tilde{\eta}, \label{eq:engine_avcpg_f} \\
    \calI_n[\tilde{T}_c\tilde{W}_c] &= \int_{T_n} T_c\tilde{W}_c,  \qquad  c = 1, \dots, C, \label{eq:engine_avcpg_g}
\end{align}
\end{subequations}
where, for each $c$, $P_c$ and $T_c$ are functions of $S_c$ and $V_c = V_p - \cos(\theta - 2\pi c/C)$ \eqref{eq:piston_volume} given by the working fluid's constitutive relations.
By \Cref{th:generic_ode_sp}, the scheme \eqref{eq:engine_avcpg} conserves $E$ and is increasing in $S$ (up to machine precision, quadrature and solver tolerances).
We observe in fact that the auxiliary variable $\tilde{\omega}$ may be eliminated from the discretisation, by reducing \eqref{eq:engine_avcpg_a} and \eqref{eq:engine_avcpg_f} to
\begin{subequations}
\begin{equation}
    \calI_n[\dot{\theta}\eta] = \int_{T_n} \omega\eta;
\end{equation}
similarly, $\tilde{P}$ may be eliminated by reducing \eqref{eq:engine_avcpg_b} and \eqref{eq:engine_avcpg_e} to
\begin{equation}
    \calI_n[\dot{\omega}\psi] = \int_{T_n} \sum_{c=1}^C P_c\sin\left(\theta - \frac{2\pi c}{C}\right)\psi.
\end{equation}
\end{subequations}

\revised{%
Observe that the cylinder temperatures $(T_c)$ in \eqref{eq:engine_avcpg_g} are positive by parametrisation:
$T_c = \exp(S_c/C_V)V_c^{-1/C_V}$ is positive for any value of $S_c$ (since $V_c \ge V_p - 1 > 0$ \eqref{eq:piston_volume}).
This guarantee does not hold for the auxiliary temperatures $(\tilde{T}_c)$.
As they enter the discretisation \eqref{eq:engine_avcpg} through $T_0/\tilde{T}_c$, their vanishing would render $\tilde{D}$ singular, something we are unable to preclude \emph{a priori}.
The auxiliary temperatures in each of our tests below remain well bounded away from zero; however, were this to cause an issue in a practical application, the timestep could simply be reduced adaptively.}

\subsubsection{Numerical test}

To demonstrate the discretisation, we close the system \eqref{eq:engine_system} and its discretisation \eqref{eq:engine_avcpg} with a fixed equation of state, that of an ideal fluid with heat capacity $C_V$ and adiabatic index $\gamma = 1 + 1/C_V$:
\begin{equation}
    P_c = \exp\left(\frac{S_c}{C_V}\right) V_c^{-\gamma},  \qquad
    T_c = P_c V_c,  \qquad
    U_c = C_V T_c,  \qquad
    U_0 = T_0 S_0.
\end{equation}
For our experiments, we consider $C_V = 5/2$.
We further take $C = 6$ cylinders with average piston volume $V_p = 1 + 2^{-4}$, and assume the system is initially in thermal equilibrium with $T_c(0) = T_0 = 1$ for each $c$.
The initial angular displacement and external entropy are both set to $0$, while the initial rate of change of the engine phase is chosen as $\omega(0) = 8$.

\begin{figure}[!ht]
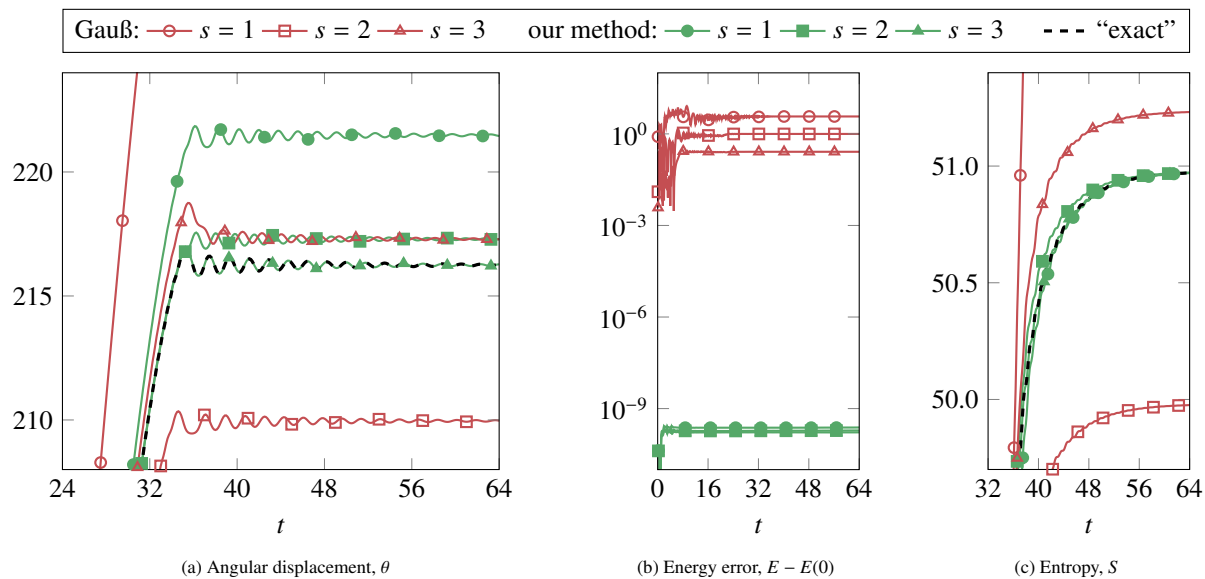

    \captionsetup[subfigure]{justification = centering}
    \centering

    \begin{subfigure}{0.45\textwidth}
        \centering



        \caption{Angular displacement, $\theta$}
        \label{fig:engine_short_angle}
    \end{subfigure}
    \begin{subfigure}{0.26\textwidth}
        \centering

        %

        \caption{Energy error, $E - E(0)$}
    \end{subfigure}
    \begin{subfigure}{0.26\textwidth}
        \centering

        %

        \caption{Entropy, $S$}
        \label{fig:engine_short_entropy}
    \end{subfigure}

    \caption{Evolution in $\theta$, $E$, $S$ for the combustion engine at $\Delta t_n = 2^{-4}$.}
    \label{fig:engine_short}
\end{figure}

\begin{figure}[!ht]
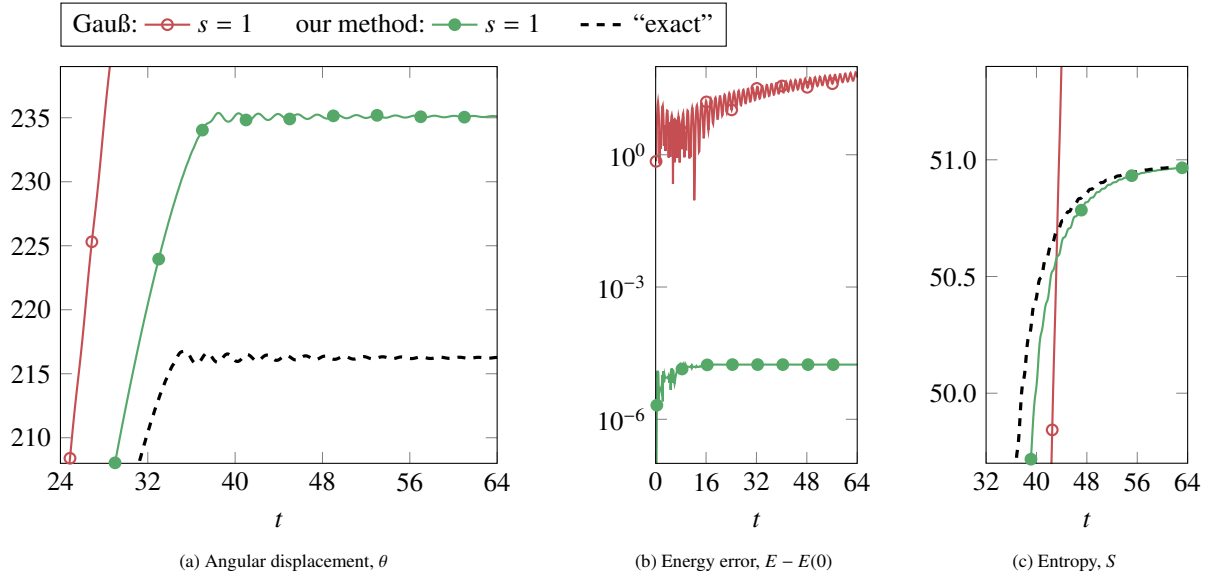

    \captionsetup[subfigure]{justification = centering}
    \centering

    \begin{subfigure}{0.45\textwidth}
        \centering



        \caption{Angular displacement, $\theta$}
        \label{fig:engine_long_angle}
    \end{subfigure}
    \begin{subfigure}{0.26\textwidth}
        \centering

        %

        \caption{Energy error, $E - E(0)$}
    \end{subfigure}
    \begin{subfigure}{0.26\textwidth}
        \centering

        %

        \caption{Entropy, $S$}
    \end{subfigure}

    \caption{Evolution in $\theta$, $E$, $S$ for the combustion engine at $\Delta t_n = 2^{-3}$.}
    \label{fig:engine_long}
\end{figure}

Figures~\ref{fig:engine_short} \&~\ref{fig:engine_long} compare the evolution of $\theta$, $E$ and $S$ in numerical simulations of this model engine \eqref{eq:engine_system} with both our proposed energy- and entropy-stable integrator \eqref{eq:engine_avcpg} taking $\calI_n$ to be an $s$-stage Gau\ss{}--Legendre method, and the comparable (symplectic) Gau\ss{} method.
Each simulation is run up to a final time $t = 2^6$, shown alongside a finely resolved ``exact'' trajectory computed using our integrator with $s=4$ stages and a small timestep $\Delta t_n = 2^{-5}$.

\Cref{fig:engine_short} uses a timestep $\Delta t_n = 2^{-4}$, with stages $s \in \{1, 2, 3\}$.
We see from \Cref{fig:engine_short_angle} that our method \eqref{eq:engine_avcpg} at lowest order $s=1$ has a smaller final error in the angular displacement $\theta$ than both the 1- and 2-stage Gau\ss{} methods;
at order $s=3$, there is no visible difference between the results from our scheme and the ``exact'' trajectory.
From \Cref{fig:engine_short_entropy}, we see that our lowest-order ($s=1$) simulation more accurately reproduces the evolution of the entropy $S$ than the highest-order ($s=3$) simulation using the Gau\ss{} method.
Despite the symplectic nature of the Gau\ss{} methods, each of the numerical trajectories exhibits a notable final error in the energy $E$, approximately $3.7$, $1.0$ and $0.3$ for stages $s=1$, $2$ and $3$ respectively (compare with the initial value of $E(0) \approx 67.8$).

\Cref{fig:engine_long} uses a longer timestep $\Delta t_n = 2^{-3}$.
We consider only the lowest order $s=1$ as the Gau\ss{} method fails with stages $s > 1$.
This represents no major cause for concern, as we have no reason to expect a numerical solution to exist for such implicit multistage methods at such long timesteps.
\revised{%
Under an $\ell^2$ line search (in place of the default Newton solver used throughout), our integrator \eqref{eq:engine_avcpg} \emph{does} converge at this timestep with $s=2$ (though still not with $s=3$).
Over a shorter duration $t \in [0, 16]$, our integrator also remains solvable at $s=1$ up to $\Delta t_n = 2^{-2}$, whereas the Gau\ss{} method does not\footnote{
    \revised{%
    As noted in \Cref{rem:stiff_problems}, over increasingly long timesteps an L-stable discretisation may be preferred.
    For faster heat exchange, the restriction tightens accordingly.
    If the friction matrix $D$ is scaled by a factor $\kappa$ (which preserves the energy-conserving and entropy-generating properties of the scheme), the thermal timescale comes to dominate the mechanical one;
    for $\kappa$ ranging from $2^4$ to $2^{10}$, we find the nonlinear solve succeeds when $\kappa\Delta t_n \le 2^2$ and fails when $\kappa\Delta t_n \ge 2^3$.
    }
}.}

Unlike the numerical results from our method, the angular displacement $\theta$, energy $E$ and entropy $S$ in the implicit midpoint simulation all climb steadily throughout the duration, as the symplectic scheme fails to capture the dissipative effects of the heat exchange.
The angle $\theta$ reaches a final value of approximately $537$ by the end of the simulation;
the energy $E$ reaches approximately $125.4$ (compare again with $E(0) \approx 67.8$), while the entropy $S$ reaches a final value of approximately $73.8$ (compare with the initial value $S(0) \approx 2.3$ and ``exact'' final value of approximately $51.0$).

\revised{%
In practice we observe that each auxiliary temperature $\tilde{T}_c$ remains positive throughout each of the simulations reported above, ensuring $\tilde{D}$ remains non-singular.
The smallest value across them, approximately $0.558$ against an external temperature $T_0 = 1$, is attained during the first compression of the 2-stage method at $\Delta t_n = 2^{-4}$.}

\section{GENERIC PDEs}\label{sec:generic_pdes}

Having developed the structure-preserving discretisations of finite-dimensional GENERIC systems, we now extend the same thermodynamically consistent framework to the infinite-dimensional setting, with the key ideas carrying over upon replacement of the finite-dimensional Euclidean space $\bbR^d$ with a Banach space $U$.
The GENERIC PDE is most conveniently stated in a variational form analogous to \eqref{eq:generic_ode}:
find $u \in C^1(\bbR_+; U)$ satisfying known initial data, such that
\begin{equation}\label{eq:generic_pde}
    M(u; \dot{u}, v)  =  B(u, w_E(u); v) + D(u, w_S(u); v)
\end{equation}
for all test functions $v \in U$ and all times $t \in \bbR_+$. Here $M : U \times U \times U \to \bbR$ is bilinear in its final two arguments. The operators $B, D : U \times U \times U \to \bbR$ are each linear in the final argument $v$\footnote{This is the significance of the semicolon. The GENERIC framework typically requires that $B, D : U \times U \times U \to \bbR$ each be bilinear in their final two arguments; however, this condition is not necessary for the construction of our energy- and entropy-preserving scheme.}; $B(u, \cdot; \cdot)$ and $D(u, \cdot; \cdot)$ are, respectively, skew-symmetric and positive-semidefinite (the Poisson and friction operators). The functionals $E, S : U \to \bbR$ denote the energy and entropy.
The functionals $w_E, w_S : U \to U$ are such that $M(u; \cdot, w_E(u)) = E'(u; \cdot)$, $M(u; \cdot, w_S(u)) = S'(u; \cdot)$, where $E'$, $S'$ are the Fr\'echet derivatives of $E$, $S$ respectively.
As in \eqref{eq:generic_ode_compatibility}, thermodynamic consistency imposes the conditions
\begin{align}\label{eq:generic_pde_compatibility}
    B(u, \cdot; w_S(u))  =  0,  \qquad
    D(u, \cdot; w_E(u))  =  0.
\end{align}
We see these conditions ensure energy conservation ($\tfrac{\rmd}{\rmd t}E(u) = 0$) and entropy generation ($\tfrac{\rmd}{\rmd t}S(u) \ge 0$) along exact trajectories by considering $v = w_E(u)$ and $v = w_S(u)$ in \eqref{eq:generic_pde};
accordingly, $w_E(u)$ and $w_S(u)$ are the associated test functions for energy $E$ and entropy $S$, which we shall introduce into our discretisation through auxiliary variables, i.e.~projections $\tilde{w}_E$ and $\tilde{w}_S$ respectively onto a discrete test space.
As in \Cref{ass:generic_matrices}, for our construction of an energy- and entropy-stable integrator we assume access to certain extensions of $B$ and $D$.

\begin{assumption}[Characterisation of GENERIC operator compatibility]\label{ass:generic_operators}
    Assume the existence of $\tilde{B}, \tilde{D} : U \times U \times U \times U \to \bbR$, linear in each of their final arguments, such that the following hold:
    \begin{subequations}
    \begin{enumerate}
        \item Consistency with $B, D$:
        for all $u$,
        \begin{align}\label{eq:generic_operator_coincidence}
            \tilde{B}(u, w_S(u), \cdot; \cdot) = B(u, \cdot; \cdot),  \qquad
            \tilde{D}(u, w_E(u), \cdot; \cdot) = D(u, \cdot; \cdot).
        \end{align}
        \item Skew-symmetry and positive semidefiniteness:
        for all $u, \tilde{w}_E, \tilde{w}_S \in U$, $\tilde{B}(u, \tilde{w}_S, \cdot; \cdot)$ is skew-symmetric, and $\tilde{D}(u, \tilde{w}_E, \cdot; \cdot)$ is positive-semidefinite.
        \item Preservation of compatibility \eqref{eq:generic_pde_compatibility}:
        for all $u, \tilde{w}_E, \tilde{w}_S \in U$,
        \begin{align}\label{eq:generic_operator_compatibility}
            \tilde{B}(u, \tilde{w}_S, \cdot; \tilde{w}_S)  =  0,  \qquad
            \tilde{D}(u, \tilde{w}_E, \cdot; \tilde{w}_E)  =  0.
        \end{align}
    \end{enumerate}
    \end{subequations}
\end{assumption}

\begin{remark}\label{rem:pde_scope}
    \revised{%
    The considerations of \Cref{rem:ode_scope}, regarding when such compatibility conditions may be satisfied, apply again here.
    There is an additional subtlety in the PDE setting, in that the continuous compatibility conditions may hold through integration by parts, as in the Boltzmann example below (\Cref{sec:boltzmann}), and so their discrete counterparts can further demand care with boundary terms, weak formulations, and choices of discrete spaces.
    }
\end{remark}

With \Cref{ass:generic_operators} in place, we follow a similar strategy to each of the ODE cases above to construct an energy- and entropy-stable integrator for \eqref{eq:generic_pde}, by (i)~introducing auxiliary variables $\tilde{w}_E$ and $\tilde{w}_S$ approximating $w_E(u)$ and $w_S(u)$, and (ii)~introducing these auxiliary variables into the right-hand side of \eqref{eq:generic_pde}.
Discretising in space, let $\bbU$ be a conforming finite-dimensional subspace of $U$, e.g.~a finite element space.
Defining the spacetime trial space for a given timestep
\begin{equation}\label{eq:solution_space}
    \bbX_n  \coloneqq  \left\{ \, u \in \bbP_s(T_n; \bbU) \mid u(t_n) \text{ satisfies known initial data} \, \right\},
\end{equation}
our structure-preserving scheme is:
find $(u, (\tilde{w}_E, \tilde{w}_S)) \in \bbX_n \times \dot{\bbX}_n^2$ such that for all $(v, (v_E, v_S)) \in \dot{\bbX}_n \times \dot{\bbX}_n^2$,
\begin{subequations}\label{eq:generic_pde_avcpg}
\begin{align}
    \calI_n[M(u; \dot{u}, v)]  &=  \calI_n[\tilde{B}(u, \tilde{w}_S, \tilde{w}_E; v) + \tilde{D}(u, \tilde{w}_E, \tilde{w}_S; v)],  \label{eq:generic_pde_avcpg_a}  \\
    \calI_n[M(u; v_E, \tilde{w}_E)]  &=  \int_{T_n} E'(u; v_E),  \\
    \calI_n[M(u; v_S, \tilde{w}_S)]  &=  \int_{T_n} S'(u; v_S).
\end{align}
\end{subequations}
By \eqref{eq:generic_operator_coincidence} we see \eqref{eq:generic_pde_avcpg_a} identifies with the original weak formulation \eqref{eq:generic_pde} when $(\tilde{w}_E, \tilde{w}_S) = (w_E(u), w_S(u))$.

\begin{theorem}[Energy and entropy stability of \eqref{eq:generic_pde_avcpg}]\label{th:generic_pde_sp}
    Assuming solutions exist, the integrator \eqref{eq:generic_pde_avcpg} is energy- and entropy-stable, with $E(u(t_{n+1})) = E(u(t_n))$ and $S(u(t_{n+1})) \ge S(u(t_n))$.
\end{theorem}

\begin{proof}
    By considering respectively $v_E = \dot{u}$, $v_S = \dot{u}$ and $v = \tilde{w}_E$, $v = \tilde{w}_S$ in \eqref{eq:generic_pde_avcpg},
    \begin{subequations}
    \begin{align}
        E(u(t_{n+1})) &- E(u(t_n))  &
            S(u(t_{n+1})) &- S(u(t_n))  \notag  \\
        &= \int_{T_n} \dot{E}  &
            &= \int_{T_n} \dot{S}  \\
        &= \int_{T_n} E'(u; \dot{u})  &
            &= \int_{T_n} S'(u; \dot{u})  \\
        &= \calI_n\big[M(u; \dot{u}, \tilde{w}_E)\big]  &
            &= \calI_n\big[M(u; \dot{u}, \tilde{w}_S)\big]  \\
        &= \calI_n\left[\begin{aligned}
            &\tilde{B}(u, \tilde{w}_S, \tilde{w}_E; \tilde{w}_E)  \\
            &\qquad\qquad+ \tilde{D}(u, \tilde{w}_E, \tilde{w}_S; \tilde{w}_E)
        \end{aligned}\right]  &
            &= \calI_n\left[\begin{aligned}
                &\tilde{B}(u, \tilde{w}_S, \tilde{w}_E; \tilde{w}_S)  \\
                &\qquad\qquad+ \tilde{D}(u, \tilde{w}_E, \tilde{w}_S; \tilde{w}_S)
            \end{aligned}\right]  \\
        &= 0,  &
            &\ge  0,
    \end{align}
    \end{subequations}
    where the final equality and inequality hold by \Cref{ass:generic_operators}.
\end{proof}

\begin{remark}\label{rem:generic_pde_cost}
    Unlike in the finite-dimensional cases above, the auxiliary variables $\tilde{w}_E$, $\tilde{w}_S$ generally cannot be eliminated.
    The additional cost incurred in solving for them is nontrivial.
\end{remark}

\revised{%
As in the ODE setting, these stability properties are conditional on the existence of discrete solutions.
A full existence and convergence theory in the present PDE setting would additionally require an analysis of the spatial discretisation, which we have deliberately left unspecified (assuming only a conforming finite-dimensional subspace $\bbU \subset U$).
Such an analysis is necessarily problem-dependent, and lies beyond the scope of the present work.
We therefore restrict our attention to the structure-preserving properties, which hold irrespective of the spatial discretisation.
}

\subsection{Boltzmann equation}\label{sec:boltzmann}

The Boltzmann equation is a key example of a GENERIC PDE.
\revised{%
We include it as a formal demonstration of the abstract construction of \Cref{th:generic_pde_sp} applied to a challenging and physically important example.
A fully practical numerical discretisation of the Boltzmann equation lies beyond the scope of the present work.}

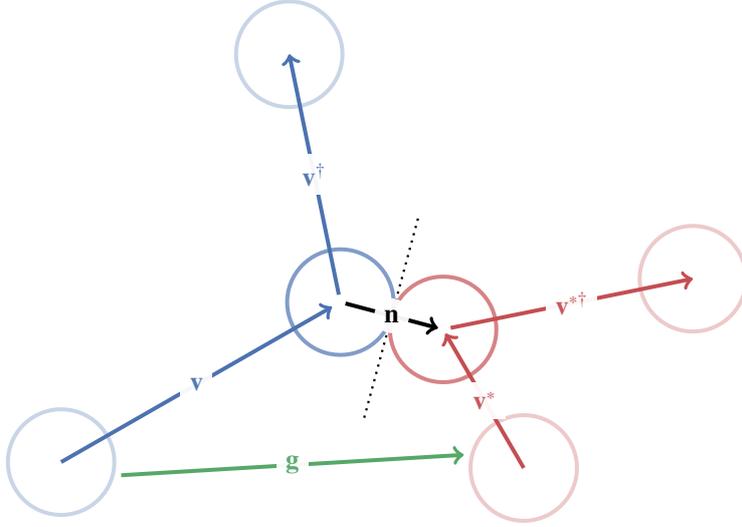
\begin{figure}[!ht]
    \centering

    \begin{tikzpicture}[line width=1.5pt]
        \begin{scope}[rotate=75]
        \def\scale{0.3}
        \def\r{0.7}

        \coordinate (A) at (0,\r);
        \coordinate (B) at (0,-\r);

        \coordinate (v) at (10,-10);
        \coordinate (v*) at (5,5);
        \coordinate (vdag) at (10,5);
        \coordinate (v*dag) at (5,-10);

        \coordinate (A-pre) at ($(A)-\scale*(v)$);
        \coordinate (B-pre) at ($(B)-\scale*(v*)$);
        \coordinate (A-post) at ($(A)+\scale*(vdag)$);
        \coordinate (B-post) at ($(B)+\scale*(v*dag)$);

        \draw[seabornblue!70] (A) circle [radius=\r];
        \draw[seabornred!70] (B) circle [radius=\r];
        \draw[seabornblue!30] (A-pre) circle [radius=\r];
        \draw[seabornred!30] (B-pre) circle [radius=\r];
        \draw[seabornblue!30] (A-post) circle [radius=\r];
        \draw[seabornred!30] (B-post) circle [radius=\r];

        \coordinate (g) at ($(v) - (v*)$);
        \coordinate (g-mid) at ($(A-pre)!0.5!(B-pre)$);
        \draw[seaborngreen, ->] ($(g-mid) - \scale/2.1*(g)$) -- ($(g-mid) + \scale/2.1*(g)$)
            node[midway, circle, inner sep=1pt, fill=white, fill opacity=0.95] {$\bfg$};

        \draw[->] ($(A)!0.05!(B)$) -- ($(A)!0.95!(B)$);
        \draw[dotted, line width=0.85pt] ($(-2*\r,0)$) -- ($(2*\r,0)$);
        \draw[draw=none, fill=white, fill opacity=0.95] (0,0) circle [radius=0.23];
        \node at (0,0) {$\bfn$};

        \draw[seabornblue, ->] (A-pre) -- ($(A)!0.1*\scale!(A-pre)$)
            node[midway, circle, inner sep=1pt, fill=white, fill opacity=0.95] {$\bfv$};
        \draw[seabornred, ->] (B-pre) -- ($(B)!0.1*\scale!(B-pre)$)
            node[midway, circle, inner sep=1pt, fill=white, fill opacity=0.95] {$\bfv^*$};
        \draw[seabornblue, ->] ($(A)!0.1*\scale!(A-post)$) -- (A-post)
            node[midway, circle, inner sep=1pt, fill=white, fill opacity=0.95] {$\bfv^\dag$};
        \draw[seabornred, ->] ($(B)!0.1*\scale!(B-post)$) -- (B-post)
            node[midway, circle, inner sep=1pt, fill=white, fill opacity=0.95] {$\bfv^{*\dag}$};
        \end{scope}
    \end{tikzpicture}

    \caption{Illustration of the relationship between the velocities $\bfv$, $\bfv^*$, $\bfv^\dag$, $\bfv^{*\dag}$, and collision normal $\bfn$.}
    \label{fig:collision}
\end{figure}

We consider the nondimensionalised Boltzmann equation in $d \ge 1$ dimensions,
\begin{equation}\label{eq:boltzmann}
    \dot{f}  =  - \bfv\cdot\nabla_\bfx f + \nabla_\bfx\phi\cdot\nabla_\bfv f + \frac{1}{\Kn}C(f).
\end{equation}
Here $f(\bfx, \bfv, t) \in \bbR$ represents the particle density function at position $\bfx \in \bbR^d$ and velocity $\bfv \in \bbR^d$, $\nabla_\bfx$ and $\nabla_\bfv$ denote the partial derivatives with respect to $\bfx$ and $\bfv$, $\phi(\bfx) \in \bbR$ represents a potential energy density, and $\Kn$ is the Knudsen number.
The term $C$ denotes the Boltzmann collision operator, defined by
\begin{equation}
    C(f)  \coloneqq  \int_{\bfv^*, \bfn \in S^{d-1}}\beta(\theta, \|\bfg\|)\,(f^\dag f^{*\dag} - ff^*) \, \rmd\bfv^*\rmd\bfn.
\end{equation}
Here $\bfg \coloneqq \bfv - \bfv^*$ is the relative velocity, $\bfn$ is the collision normal (on the unit $(d-1)$-sphere $S^{d-1} \subset \bbR^d$), $\theta = 2\angle(\bfg, \bfn)$ is the scattering angle, and $\beta(\theta, \|\bfg\|) \ge 0$ is the collision kernel, where $\|\cdot\|$ denotes the $\ell^2$ norm;
the fields $f^*$, $f^\dag$, $f^{*\dag}$ are shorthand for
\begin{subequations}
\begin{equation}\label{eq:boltzmann_f_shorthand}
    f^*        =  f|_{\bfv = \bfv^*},   \qquad
    f^\dag     =  f|_{\bfv = \bfv^\dag},    \qquad
    f^{*\dag}  =  f|_{\bfv = \bfv^{*\dag}},
\end{equation}
where in turn $\bfv^\dag, \bfv^{*\dag}$ are the unique post-collision velocities (see \Cref{fig:collision}) satisfying
\begin{equation}\label{eq:boltzmann_collision_vectors}
    \bfv + \bfv^*  =  \bfv^\dag + \bfv^{*\dag},  \qquad
    \frac{1}{2}\|\bfv\|^2 + \frac{1}{2}\|\bfv^*\|^2  =  \frac{1}{2}\|\bfv^\dag\|^2 + \frac{1}{2}\|\bfv^{*\dag}\|^2,  \qquad
    \bfn  =  \frac{\bfv - \bfv^\dag}{\|\bfv - \bfv^\dag\|}.
\end{equation}
\end{subequations}
We assume periodic boundary conditions in $\bfx$, and a vanishing asymptotic boundary condition in $\bfv$ of $f \to 0$ as $\|\bfv\| \to \infty$.
The Boltzmann equation \eqref{eq:boltzmann} has a conserved energy $E$ and generated entropy $S$,
\begin{equation}\label{eq:boltzmann_qois}
    E  \coloneqq  \int_{\bfx, \bfv}\left(\frac{1}{2}\|\bfv\|^2 + \phi\right)f,  \qquad
    S  \coloneqq  \int_{\bfx, \bfv}(1 - \log f)f.
\end{equation}

To handle the asymptotic boundary conditions in $\bfv$, we parametrise $f$ as
\begin{equation}\label{eq:boltzmann_f}
    f(\bfx, \bfv, t)  =  f_0(\bfv)\exp(u(\bfx, \bfv, t))
\end{equation}
for $u$ in some sufficiently regular function space $U$.
The function $f_0 > 0$ characterises the asymptotic behaviour in $\bfv$, with $f_0 \to 0$ and $u = o[\log f_0]$ as $\|\bfv\| \to \infty$.
Note then that $\dot{f} = f\dot{u}$.

To apply the scheme \eqref{eq:generic_pde_avcpg} to construct an energy- and entropy-stable integrator for \eqref{eq:boltzmann}, we recall the formulation of the Boltzmann equation in the GENERIC formalism presented by \"Ottinger \cite{Ottinger_1997}.

We first cast \eqref{eq:boltzmann} into a variational form.
By testing against $v \in U$ and performing some standard manipulations of the collision term, we arrive at the following:
find $u \in C^1(\bbR_+; U)$ satisfying known initial data, such that
\begin{equation}\label{eq:boltzmann_variational}
      \int_{\bfx, \bfv}f\dot{u}v
    = \int_{\bfx, \bfv}(\nabla_\bfx\phi\cdot\nabla_\bfv v - \bfv\cdot\nabla_\bfx v)f
    + \frac{1}{4\Kn}\int_{\bfx, \bfv, \bfv^*, \bfn}\beta\,(f^\dag f^{*\dag} - ff^*)\,(v + v^* - v^\dag - v^{*\dag})
\end{equation}
at all times $t \in \bbR_+$ and for all $v \in U$, where $v^*$, $v^\dag$, $v^{*\dag}$ are defined analogously to $f^*$, $f^\dag$, $f^{*\dag}$ as in \eqref{eq:boltzmann_f_shorthand}.
This induces the choice of the left-hand side operator $M : U \times U \times U \to \bbR$,
\begin{equation}
    M(u; w, v)  \coloneqq  \int_{\bfx, \bfv} fwv  =  \int_{\bfx, \bfv} f_0\exp(u)wv.
\end{equation}

Now, consider the energy $E$ and entropy $S$ \eqref{eq:boltzmann_qois} as functionals of $u$.
These have Fr\'echet derivatives
\begin{equation}
    E'(u; \delta u)  =  \int_{\bfx, \bfv}\left(\frac{1}{2}\|\bfv\|^2 + \phi\right)f \delta u,  \qquad
    S'(u; \delta u)  =  - \int_{\bfx, \bfv} f \log f \delta u,
\end{equation}
where again $f$ is defined in terms of $u$ by \eqref{eq:boltzmann_f}.
Seeking $w_E(u)$, $w_S(u)$ such that $M(u; \cdot, w_E(u)) = E'(u; \cdot)$, $M(u; \cdot, w_S(u)) = S'(u; \cdot)$, the solution is immediate:
\begin{align}
    w_E(u)  =  \frac{1}{2}\|\bfv\|^2 + \phi,   \qquad
    w_S(u)  =  - \log f.
\end{align}
Up to a constant factor in $w_S$, these align with \cite[eqs.~(6,~7)]{Ottinger_1997}.
Define the Poisson and friction operators $B, D : U \times U \times U \to \bbR$,
\begin{subequations}
\begin{align}
    B(u, \tilde{w}_E; v)
        &\coloneqq  \int_{\bfx, \bfv}(\nabla_\bfx \tilde{w}_E\cdot\nabla_\bfv v - \nabla_\bfv \tilde{w}_E\cdot\nabla_\bfx v)f,  \label{eq:boltzmann_poisson}  \\
    D(u, \tilde{w}_S; v)
        &\coloneqq  \frac{1}{4\Kn}\int_{\bfx, \bfv, \bfv^*, \bfn}\beta\left(\exp(- \tilde{w}_S^\dag - \tilde{w}_S^{*\dag}) - \exp(- \tilde{w}_S - \tilde{w}_S^*)\right)(v + v^* - v^\dag - v^{*\dag}),  \label{eq:boltzmann_friction}
\end{align}
\end{subequations}
where $\tilde{w}_S^*$, $\tilde{w}_S^\dag$, $\tilde{w}_S^{*\dag}$ are again defined analogously to $f^*$, $f^\dag$, $f^{*\dag}$ \eqref{eq:boltzmann_f_shorthand} and $v^*$, $v^\dag$, $v^{*\dag}$ \eqref{eq:boltzmann_collision_vectors}.
The Poisson operator $B$ in \eqref{eq:boltzmann_poisson} aligns immediately with the operator $L$ as defined in \cite[eq.~(8)]{Ottinger_1997}.
The skew-symmetry of $B$ is immediate, while the positive-semidefiniteness of $D$ relies on the observation that $(e^{-x} - e^{-y})(y - x) \ge 0$.
The GENERIC compatibility condition $B(u, \cdot; {w}_S) = 0$ holds after algebraic simplification, while $D(u, \cdot; {w}_E) = 0$ can be seen from the conservation of energy condition $\frac{1}{2}\|\bfv\|^2 + \frac{1}{2}\|\bfv^*\|^2  =  \frac{1}{2}\|\bfv^\dag\|^2 + \frac{1}{2}\|\bfv^{*\dag}\|^2$ \eqref{eq:boltzmann_collision_vectors}.
With $M$, $B$, $D$ as defined, the Boltzmann equation is a GENERIC PDE of the form \eqref{eq:generic_pde};
we can thus apply \eqref{eq:generic_pde_avcpg} to preserve energy and entropy stability, provided we can
define suitable $\tilde{B}, \tilde{D} : U \times U \times U \times U \to \bbR$ satisfying \Cref{ass:generic_operators}.

Take $\tilde{w}_E, \tilde{w}_S$ to be approximations to $w_E(u), w_S(u)$.
The incorporation of $\tilde{w}_S$ into $B$ is simple, with $\tilde{B}$ defined by
\begin{equation}
    \tilde{B}(u, \tilde{w}_S; \tilde{w}_E, v)  \coloneqq  \int_{\bfx, \bfv}(\nabla_\bfx \tilde{w}_E\cdot\nabla_\bfv v - \nabla_\bfv \tilde{w}_E\cdot\nabla_\bfx v)\exp(- \tilde{w}_S).
\end{equation}
In the case of $\tilde{w}_S = w_S(u) = - \log f$, we see this identifies with $B$ as $\exp(- \tilde{w}_S) = \exp(\log f) = f$;
we see $\tilde{B}$ evaluates to $0$ for $v = \tilde{w}_S$ by the substitution $\exp(- \tilde{w}_S)\nabla\tilde{w}_S = - \nabla[\exp(- \tilde{w}_S)]$ and integration by parts in $\bfx$, $\bfv$ noting the periodic boundary conditions and assuming $\tilde{w}_S \to \infty$ as $\|\bfv\| \to \infty$ such that $\exp(- \tilde{w}_S) \to 0$.

The incorporation of $\tilde{w}_E$ into $D$ is somewhat more involved, and requires rewriting \eqref{eq:boltzmann_friction}.
For concreteness we fix $d = 3$ in what follows.
Let $\Sigma \subset \bbR^3 \times (\bbR^3)^4$ denote the manifold of tuples $(\bfx, (\bfv, \bfv^*, \bfv^\dag, \bfv^{*\dag}))$ satisfying the relations \eqref{eq:boltzmann_collision_vectors} characterising all possible collisions;
the final relation in \eqref{eq:boltzmann_collision_vectors} now instead serves as a definition for $\bfn = \tfrac{\bfv - \bfv^\dag}{\|\bfv - \bfv^\dag\|}$ and consequently $\theta = 2\angle(\bfg, \bfn)$.
We endow $\Sigma$ with the metric induced from $\bbR^3 \times (\bbR^3)^4$.
The friction operator $D$ may then be written as an integral over $\Sigma$,
\begin{equation}
    D(u; \tilde{w}_S, v)
        \coloneqq  \frac{1}{4\Kn}\int_\Sigma \beta\left(\exp(- \tilde{w}_S^\dag - \tilde{w}_S^{*\dag}) - \exp(- \tilde{w}_S - \tilde{w}_S^*)\right)
        (v + v^* - v^\dag - v^{*\dag}).
\end{equation}
This formulation of the friction operator aligns with \cite[eq.~(12)]{Ottinger_1997}, up to a minor modification made by \"Ottinger to enforce linearity\footnote{This modification is unnecessary for our purposes as the structure-preserving properties hold regardless. \"Ottinger \cite{Ottinger_1997} characterises both the collision kernel $\beta$ and the manifold of admissible collisions $\Sigma$ through a single distribution function $w$ of transition probabilities supported on $\Sigma$.}.

With this reformulation, we may define an auxiliary manifold $\tilde{\Sigma} \subset \bbR^3 \times (\bbR^3)^4$ of tuples $(\bfx, (\bfv, \bfv^*, \bfv^{\tilde{\dag}}, \bfv^{*\tilde{\dag}}))$ satisfying the auxiliary relations
\begin{subequations}
\begin{equation}\label{eq:boltzmann_collision_vectors_avcpg}
    \bfv + \bfv^*  =  \bfv^{\tilde{\dag}} + \bfv^{*\tilde{\dag}},  \qquad
    \tilde{w}_E|_{\bfv = \bfv} + \tilde{w}_E|_{\bfv = \bfv^*}  =  \tilde{w}_E|_{\bfv = \bfv^{\tilde{\dag}}} + \tilde{w}_E|_{\bfv = \bfv^{*\tilde{\dag}}},
\end{equation}
which we again endow with the metric induced from $\bbR^3 \times (\bbR^3)^4$.
Similarly to $f^*, f^\dag, f^{*\dag}$ as in \eqref{eq:boltzmann_f_shorthand} we take $\psi^*, \psi^{\tilde{\dag}}, \psi^{*\tilde{\dag}}$, for an arbitrary function $\psi$ in $\bfv$, as shorthand for
\begin{equation}\label{eq:boltzmann_f_shorthand_avcpg}
    \psi^*                \coloneqq  \psi|_{\bfv = \bfv^*},  \qquad
    \psi^{\tilde{\dag}}   \coloneqq  \psi|_{\bfv = \bfv^{\tilde{\dag}}},  \qquad
    \psi^{*\tilde{\dag}}  \coloneqq  \psi|_{\bfv = \bfv^{*\tilde{\dag}}}.
\end{equation}
\end{subequations}
We may then introduce $\tilde{w}_E$ into the definition of $\tilde{D}$ implicitly through $\tilde{\Sigma}$:
\begin{equation}
    \tilde{D}(u, \tilde{w}_E; \tilde{w}_S, v)
        \coloneqq  \frac{1}{4\Kn}\int_{\tilde{\Sigma}} \beta\left(\exp(- \tilde{w}_S^{\tilde{\dag}} - \tilde{w}_S^{*\tilde{\dag}}) - \exp(- \tilde{w}_S - \tilde{w}_S^*)\right)
        (v + v^* - v^{\tilde{\dag}} - v^{*\tilde{\dag}}).
\end{equation}
Considering $\tilde{w}_E = w_E(u) = \frac{1}{2}\|\bfv\|^2 + \phi(\bfx)$, we see this identifies with $D$ as the conditions \eqref{eq:boltzmann_collision_vectors_avcpg} on $\bfv$, $\bfv^*$, $\bfv^{\tilde{\dag}}$, $\bfv^{*\tilde{\dag}}$ align with those on $\bfv$, $\bfv^*$, $\bfv^\dag$, $\bfv^{*\dag}$ \eqref{eq:boltzmann_collision_vectors};
we see $\tilde{D}$ evaluates to zero for $v = \tilde{w}_E$ as the conditions \eqref{eq:boltzmann_collision_vectors_avcpg} on $\bfv^*$, $\bfv^{\tilde{\dag}}$, $\bfv^{*\tilde{\dag}}$ imply $\tilde{w}_E + \tilde{w}_E^* - \tilde{w}_E^{\tilde{\dag}} - \tilde{w}_E^{*\tilde{\dag}} = 0$ on $\tilde{\Sigma}$.

For a finite-dimensional subspace $\bbU \subset U$ and discrete space-time function space $\bbX_n$ defined as in \eqref{eq:solution_space}\footnote{Since the Boltzmann equation has a hyperbolic nature, a non-conforming space $\bbU$ and a corresponding non-conforming modification to \eqref{eq:boltzmann_avcpg} may be more effective in practice.}, we finally derive the following energy- and entropy-stable scheme for the Boltzmann equation \eqref{eq:boltzmann}:
find $(u, (\tilde{w}_E, \tilde{w}_S)) \in \bbX_n \times \dot{\bbX}_n^2$ such that for all $(v, (v_E, v_S)) \in \dot{\bbX}_n \times \dot{\bbX}_n^2$,
\begin{subequations}\label{eq:boltzmann_avcpg}
\begin{align}
    \calI_n\left[\int_{\bfx, \bfv}f\dot{u}v\right]  &=  \calI_n\left[\int_{\bfx, \bfv}(\nabla_\bfx\tilde{w}_E\cdot\nabla_\bfv v - \nabla_\bfv\tilde{w}_E\cdot\nabla_\bfx v)\tilde{f} + \frac{1}{4\Kn}\int_{\tilde{\Sigma}} \beta\,(\tilde{f}^{\tilde{\dag}}\tilde{f}^{*\tilde{\dag}} - \tilde{f}\tilde{f}^*)\,(v + v^* - v^{\tilde{\dag}} - v^{*\tilde{\dag}})\right],  \\
    \calI_n\left[\int_{\bfx, \bfv}f\tilde{w}_Ev_E\right]  &=  \int_{T_n}\int_{\bfx, \bfv}f\left(\frac{1}{2}\|\bfv\|^2 + \phi(\bfx)\right)v_E,  \label{eq:boltzmann_avcpg_b}  \\
    \calI_n\left[\int_{\bfx, \bfv}f\tilde{w}_Sv_S\right]  &=  - \int_{T_n}\int_{\bfx, \bfv}f\log f v_S,
\end{align}
\end{subequations}
where again $f$ is defined as in \eqref{eq:boltzmann_f}, the auxiliary density function $\tilde{f}$ is shorthand for $\tilde{f} \coloneqq \exp(- \tilde{w}_S)$, the functions $\tilde{f}^*, \tilde{f}^{\tilde{\dag}}, \tilde{f}^{*\tilde{\dag}}$ and $v^*$, $v^{\tilde{\dag}}$, $v^{*\tilde{\dag}}$ are defined via \eqref{eq:boltzmann_f_shorthand_avcpg} for $\bfv^*$, $\bfv^{\tilde{\dag}}$, $\bfv^{*\tilde{\dag}}$ defined as in \eqref{eq:boltzmann_collision_vectors_avcpg}, and $\tilde{\Sigma} \subset \bbR^3 \times (\bbR^3)^4$ is the auxiliary manifold of $(\bfx, (\bfv, \bfv^*, \bfv^{\tilde{\dag}}, \bfv^{*\tilde{\dag}}))$ satisfying \eqref{eq:boltzmann_collision_vectors_avcpg}.

As in the proof of \Cref{th:generic_pde_sp}, the conservation of $E$ and dissipation of $S$ can then be shown by testing with $(v, v_E) = (\tilde{w}_E, \dot{u})$ and $(v, v_S) = (\tilde{w}_S, \dot{u})$ respectively.

\begin{remark}
    Practical implementation of the integral over the auxiliary manifold $\tilde{\Sigma}$ in $\tilde{D}$ would in general be difficult.
    However, under certain assumptions on the space $\bbU$, this problem reduces to one only as difficult as evaluating $D$ \eqref{eq:boltzmann_friction}, as we now describe.

    Suppose $\bbU$, itself a finite-dimensional space over $(\bfx, \bfv)$, is chosen to be the tensor product $\bbU_\bfx \otimes \bbU_\bfv$ of finite-dimensional spaces over $\bfx$ and $\bfv$ respectively.
    If both $\phi \in \bbU_\bfx$ and $\frac{1}{2}\|\bfv\|^2 \in \bbU_\bfv$ (e.g.~with $\bbU_\bfv$ a degree-2 finite-element space), then $w_E = \frac{1}{2}\|\bfv\|^2 + \phi$ is an element of $\bbU$, and consequently also an element of $\dot{\bbX}_n$ by time independence.
    Choosing the quadrature rule $\calI_n = \int_{T_n}$ (i.e.~such that the scheme \eqref{eq:boltzmann_avcpg} is a structure-preserving modification of a continuous Petrov--Galerkin scheme) we see that \eqref{eq:boltzmann_avcpg_b} is then solved exactly by $\tilde{w}_E = w_E \in \dot{\bbX}_n$.
    The auxiliary variable $\tilde{w}_E$ may then be eliminated, and $\tilde{D}$ coincides with $D$.
\end{remark}

\subsection{Benjamin--Bona--Mahony equation}

In this section we consider the Benjamin--Bona--Mahony (BBM) equation \cite{Benjamin_Bona_Mahony_1972} for $u : \bbR_+ \times \Omega \to \bbR$ on an interval $\Omega \subset \bbR$,
\begin{equation}\label{eq:bbm_primal}
    \dot{u} - \partial_x^2 \dot{u}  =  - \partial_x u - u \partial_x u,
\end{equation}
where $\partial_x$ denotes the partial derivative with respect to the spatial coordinate $x$.
We impose periodic boundary conditions.
We shall observe that long-time (conservative) simulations of \eqref{eq:bbm_primal} using the simplification of \eqref{eq:generic_pde_avcpg} exhibit substantially better physical fidelity than a (symplectic) Gau\ss{} scheme of the same order with the same spatial discretisation.

The BBM equation is a Hamiltonian system with a conserved energy (denoted here by $H$ in place of $E$)
\begin{equation}
    H(u)  \coloneqq  \int_\Omega \frac{1}{2}u^2 + \frac{1}{6}u^3,
\end{equation}
and has exactly two other independent invariants:
\begin{align}\label{eq:bbm_extra_invariants}
    I_1(u) \coloneqq \int_{\Omega} u, \quad\quad
    I_2(u) \coloneqq \int_{\Omega} u^2 + u_x^2.
\end{align}
To apply the scheme \eqref{eq:generic_pde_avcpg} to construct an energy-stable integrator for \eqref{eq:bbm_primal}, we must first show \eqref{eq:bbm_primal} may be written in the form \eqref{eq:generic_pde}.

A typical semidiscrete variational form of \eqref{eq:bbm_primal} can be found by testing in the $L^2$ inner product against some test function $v$, and applying integration by parts:
find $u \in C^1(\bbR_+; U)$ satisfying known initial data such that
\begin{equation}\label{eq:bbm_variational}
    (\dot{u} - \partial_x^2 \dot{u}, v) = - \, (\partial_x u + u \partial_x u, v)
        \implies  (\dot{u}, v) + (\partial_x \dot{u}, \partial_x v) = - \, (\partial_x [u + \frac{1}{2}u^2], v)
        \implies  (\dot{u}, v)_{H^1} = (u + \frac{1}{2}u^2, \partial_x v)
\end{equation}
at all times $t \in \bbR_+$ and for all $v \in U$.
We thus take $M$ to be the $H^1$ inner product of the final two arguments.
The associated test function $w_H(u)$ for $H$ must then satisfy
\begin{equation}
    (\cdot, w_H(u))_{H^1}  =  H'(u; \cdot)
        \implies  (\cdot, w_H(u) - \partial_x^2 w_H(u))  =  (u + \frac{1}{2}u^2, \cdot),
\end{equation}
implying $w_H(u)$ is implicitly defined in strong form to satisfy
\begin{equation}
    w_H(u) - \partial_x^2 w_H(u)  =  u + \frac{1}{2}u^2.
\end{equation}
We may therefore write our semidiscrete variational form \eqref{eq:bbm_variational} in terms of $w_H(u)$ as
\begin{subequations}
\begin{align}
    (\dot{u}, v)_{H^1}
        &= (w_H(u) - \partial_x^2 w_H(u), \partial_x v)  \\
        &= \frac{1}{2}[(w_H(u), \partial_x v) + (\partial_x w_H(u), \partial_x^2 v) - (\partial_x w_H(u), v) - (\partial_x^2 w_H(u), \partial_x v)],  \\
        &= \frac{1}{2}[(w_H(u), \partial_x v)_{H^1} - (\partial_x w_H(u), v)_{H^1}].
\end{align}
\end{subequations}
If we define the skew-symmetric operator $B(w, v) \coloneqq \frac{1}{2}[(w, \partial_x v)_{H^1} \! - (\partial_x w, v)_{H^1} \!]$ (independent of $u$) and set $D = 0$, the formulation aligns with the GENERIC format \eqref{eq:generic_pde}.

For a chosen (twice weakly-differentiable) discrete space $\bbU \subset H^2(\Omega)$\footnote{With appropriate handling of arising facet terms, it is similarly possible to derive an energy-stable integrator that only requires the lower regularity $\bbU \subset H^1(\Omega)$.}, we may therefore use our general structure-preserving integrator \eqref{eq:generic_pde_avcpg} to derive the following energy-conserving scheme:
find $(u, \tilde{w}_H) \in \bbX_n \times \dot{\bbX}_n$ such that
\begin{subequations}\label{eq:bbm_avcpg}
\begin{align}
    \calI_n[(\dot{u}, v)_{H^1}]        &=  \frac{1}{2}\calI_n[(\tilde{w}_H, \partial_x v)_{H^1} \! - (\partial_x \tilde{w}_H, v)_{H^1} \!],  \\
    \calI_n[(v_H, \tilde{w}_H)_{H^1}]  &=  \int_{T_n}(u + \frac{1}{2}u^2, v_H),
\end{align}
\end{subequations}
for all $(v, v_H) \in \dot{\bbX}_n \times \dot{\bbX}_n$.
Again, taking $(v, v_H) = (\tilde{w}_H, \dot{u})$ confirms that \eqref{eq:bbm_avcpg} conserves $H$.

\subsubsection{Numerical test}\label{sec:soliton}

We simulate the very-long-time behaviour of a soliton to verify and motivate these conservation properties. We consider the domain $\Omega = (-50, 50)$.
Up to projection, the following ICs form a soliton\footnote{
    The appropriate notion of a nonlinear wave under periodic boundary conditions is not a soliton, but a cnoidal wave (see Ablowitz \& Segur \cite[Sec.~2.3]{Ablowitz_Segur_1981}).
    The values of the ICs at the boundary $x = \pm 50$ are, however, approximately $2 \times 10^{-13}$, implying that this distinction on this domain is negligible, especially after projection into the discrete space $\bbU$.
} of speed $\frac{1 + \sqrt{5}}{2}$:
\begin{equation}
    u(0)  =  \frac{3\sqrt{5} - 3}{2}\sech\left(\frac{\sqrt{5} - 1}{4}x\right)^2\!,
\end{equation}
where $\sech$ is the hyperbolic secant function.
Over an interval mesh of uniform mesh width $2$, we take $\bbU$ to be the (degree-$3$) Hermite space (see Ern \& Guermond \cite[Chap.~5]{Ern_Guermond_2021a});
in time, we take a uniform timestep $\Delta t_n = 1$.

\begin{remark}\label{rem:soliton_projection}
    \revised{%
    On the continuous level, a BBM soliton of speed $c$ satisfies $u + \frac{1}{2}u^2 = c(u - \partial_x^2 u)$, or equivalently $2 H'(u; \cdot) = c I_2'(u; \cdot)$.
    A soliton is thus a critical point of $H$ on $\hat{S} \coloneqq \{v : I_1(v) = I_1(u), \, I_2(v) = I_2(u)\}$.
    In fact, it is a constrained maximum:
    the second variation of $2 H - c I_2$ on the tangent space to $\hat{S}$ is negative definite, bar the single null direction $\partial_x u$ generated by translation, along which all three invariants are constant\footnote{
        \revised{This constrained variational characterisation is classical, and is precisely what underlies the stability of solitary waves \cite{Benjamin_1972, Grillakis_Shatah_Strauss_1987}.}
    }.
    The fully invariant manifold $S \coloneqq \{v : H(v) = H(u), \, I_1(v) = I_1(u), \, I_2(v) = I_2(u)\} \subset \hat{S}$ is therefore not of \emph{codimension} three, but of \emph{dimension} one (containing translations of the soliton only).
    Moreover, since the soliton is a constrained maximum on $\hat{S}$, initial data $u$ that is \emph{close} to being a soliton (e.g.~from spatial discretisation error) will be \emph{close} to attaining this constrained maximum;
    the fully invariant manifold $S$ will therefore be \emph{close} to one-dimensional.
    It is therefore likely that the intersection of $S$ with the discrete solution space $\bbU$ is very small, representing the difficulty of approximating translations of the discrete soliton within $\bbU$ without sacrificing at least one of the conservation laws.
    }
\end{remark}

Under these conditions, we compare the results from \revised{(i)}~a 2-stage (symplectic) Gau\ss{} method as applied to \eqref{eq:bbm_variational}\revised{, (ii)~the same 2-stage Gau\ss{} method using a $H^1$ projection\footnote{
    \revised{Following on from \Cref{rem:soliton_projection}, it is reasonable to expect that projecting onto the discrete invariant manifold $S \cap \bbU$ will be delicate.
    Indeed, we find the Newton iterations for this exact projection to fail to converge.
    We therefore relax the condition, penalising deviations in $H(u)$, $I_1(u)$ and $I_2(u)$ as heavily as possible rather than imposing conservation as a constraint. This enforces conservation up to a small tolerance on the order of $10^{-6}$.}
} to enforce conservation in $H(u)$, $I_1(u)$ and $I_2(u)$ \cite[Sec.~IV.4]{Hairer_Lubich_Wanner_2006}, and (iii)~the scheme} \eqref{eq:bbm_avcpg} with $\calI_n$ the exact integral\footnote{
    We are able to compute this exactly, as all terms in the discretisation \eqref{eq:bbm_avcpg} are polynomial. 
} and $s = 2$.

\Cref{fig:bbm_energies} shows the evolution of the energy $H(u)$ under each scheme.
Artificial dissipation of energy by the Gau\ss{} method causes the value to decrease from its initial value of around $11.1$ to around $6.2$ at the final time $t = 2 \cdot 10^4$.

\begin{figure}[!ht]
    \centering



    \caption{Evolution of the energy $H(u)$ when solving the BBM equations with a Gau\ss{} method\revised{, the same method with a $H^1$ projection for invariant conservation,} and our proposed scheme \eqref{eq:bbm_avcpg}. \revised{The projection curve is indistinguishable from that of our proposed scheme at this scale, the two differing by at most $3 \times 10^{-5}$.}}

    \label{fig:bbm_energies}
\end{figure}

\Cref{fig:bbm_plots} shows $u$ approximated with all three schemes at various times along the simulation.
The dissipation in $H(u)$ with the Gau\ss{} method\footnote{
    \revised{Note that, since the Gau\ss{} method conserves $I_1$ and $I_2$, the discrete solution is restricted to $\hat{S} \coloneqq \{v : I_1(v) = I_1(u(0)), \, I_2(v) = I_2(u(0))\}$.}
} correlates with a reduction in the amplitude of $u$, causing the speed of the soliton to decrease.
At $t = 2 \times 10^4$, the soliton in the Gau\ss{} simulation has speed approximately $1.454$;
compare with the exact value of approximately $1.618$, \deleted{and} that of the numerical solution from our proposed scheme of approximately $1.617$\revised{, and that from the projection method of approximately $1.616$}.

\begin{figure}[!ht]
    \captionsetup[subfigure]{justification = centering}
    \centering

    \begin{subfigure}{\textwidth}
        \centering



        \caption{$t = 0$}
    \end{subfigure}
    \\[3mm]
    \begin{subfigure}{\textwidth}
        \centering

        %

        \caption{$t = 1\cdot10^4$}
    \end{subfigure}
    \\[3mm]
    \begin{subfigure}{\textwidth}
        \centering

        %

        \caption{$t = 2\cdot10^4$}
    \end{subfigure}

    \caption{Plots of $u(x)$ in the BBM simulations using a Gau\ss{} method\revised{, the same method with a $H^1$ projection for invariant conservation,} and our proposed scheme at times $t \in \{0, 1 \cdot 10^4, 2 \cdot 10^4\}$. The exact solution is included for comparison. The results with the Gau\ss{} method exhibit incorrect wave speeds and oscillations throughout the domain; the results with our proposed scheme \revised{and the projection method} exhibit no such oscillations, and wave speeds accurate to three digits \revised{(albeit more accurate for our proposed scheme)}.}

    \label{fig:bbm_plots}
\end{figure}

\begin{figure}[!ht]
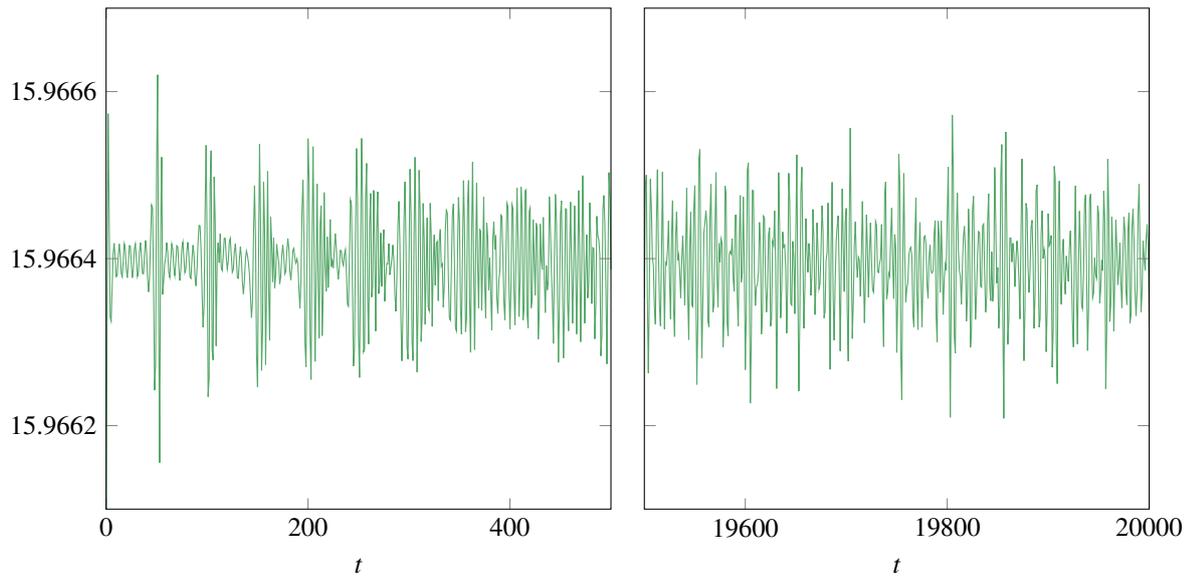

    \centering

    \begin{subfigure}{0.5\textwidth}

    \end{subfigure}

    \caption{Evolution of the squared $H^1$ norm $I_2(u)$ when solving the BBM equations with our proposed scheme.}

    \label{fig:bbm_h1}
\end{figure}

Of note is the conservation of the further invariant $I_2(u)$ \eqref{eq:bbm_extra_invariants}, i.e.~the squared $H^1$ norm $\|u\|^2_{H^1}$.
\Cref{fig:bbm_h1} shows the evolution of $I_2(u)$ under our proposed scheme.
While the construction of the scheme \eqref{eq:bbm_avcpg} does not guarantee the discrete conservation of $I_2(u)$, we find numerically that $I_2(u(t_n))$ oscillates within the small interval $(15.9660, 15.9667)$ over the simulation duration;
this is reminiscent of the approximate conservation of energy exhibited by symplectic integrators (see \Cref{fig:kepler_invariants_im} or e.g.~\cite[Chap.~I, Fig~4.1]{Hairer_Lubich_Wanner_2006}).
The proof of this property remains an open problem\footnote{
    \revised{It is somewhat surprising, given the conclusion of \Cref{rem:soliton_projection} that the discrete invariant manifold $S \cap \bbU$ may be very small and our corresponding difficulty projecting onto $S \cap \bbU$.}
}.
Conversely, as the Gau\ss{} method conserves $I_2(u)$, reductions in the $L^2$ norm of $u$ due to the decreasing energy $H(u)$ must induce increases in its $H^1$ seminorm;
this manifests as the undesirable oscillations seen throughout the domain in \Cref{fig:bbm_plots}.
No such oscillations are observed with our proposed scheme.

\section{Conclusions}\label{sec:conclusions}

In this work we have built upon our framework \cite{Andrews_Farrell_2025a} to devise conservative numerical integrators for systems of ODEs with multiple invariants, and energy- and entropy-stable integrators for ODEs and PDEs in the GENERIC formalism.

We mention two natural questions arising from this work.
The first lies in the consideration of multiply dissipative PDEs:
is it possible to use the auxiliary variable framework to preserve arbitrarily many dissipation structures for dissipative systems of ODEs or PDEs?  
The second lies in the construction of conservative integrators for multiply conservative PDEs.
Casting such equations into a semidiscrete form (discretising in space only), one retrieves a sparsely coupled system of ODEs.
Our conservative integrator for multiply conservative ODEs \eqref{eq:integrable_avcpg} could in theory be applied to this ODE system to retrieve a conservative discretisation of the original PDE. However, there is no guarantee that the alternating form construction in \Cref{sec:conservative_odes} (\Cref{th:alternating_forms}) would not damage the sparsity pattern of the system.
Is it possible to conserve arbitrarily many invariants in multi-conservative PDE systems, without breaking the sparsity structure?

\section{Code availability}

The code used for the numerical results employed Firedrake \cite{Ham_et_al_2023}, \revised{Irksome \cite{Farrell_Kirby_MarchenaMenendez_2021, Andrews_et_al_2026},} PETSc \cite{Balay_et_al_2024}, and MUMPS \cite{Amestoy_et_al_2001}.
All code for reproducing the numerical results of this work and major Firedrake components have been archived at \cite{Andrews_Farrell_2026_Zenodo}.

\section*{Acknowledgments}

This work was funded by
the European Union [ERC, GeoFEM, 101164551],
the Engineering and Physical Sciences Research Council (EPSRC) [grant number EP/W026163/1],
the EPSRC Energy Programme [grant number EP/W006839/1],
the Science and Technology Facilities Council [grant number UKRI/ST/B000495/1],
a CASE award from the UK Atomic Energy Authority,
the Donatio Universitatis Carolinae Chair ``Mathematical modelling of multicomponent systems'',
and
the UKRI Digital Research Infrastructure Programme through the Science and Technology Facilities Council's Computational Science Centre for Research Communities (CoSeC).
The authors gratefully acknowledge the hospitality of the Erwin Schrödinger International Institute for Mathematics and Physics, Vienna, where part of this work was carried out.
Views and opinions expressed are however those of the authors only and do not necessarily reflect those of the European Union or the European Research Council.
Neither the European Union nor the granting authority can be held responsible for them.
For the purpose of open access, the authors have applied a CC BY public copyright licence to any author accepted manuscript arising from this submission.
No new data were generated or analysed during this study.

\bibliographystyle{elsarticle-num}
\bibliography{references}

@book{Ablowitz_Segur_1981,
  title     = {Solitons and the Inverse Scattering Transform},
  isbn      = {978-1-61197-088-3},
  publisher = {SIAM},
  author    = {Ablowitz, M. J. and Segur, H.},
  month     = jan,
  volume    = {4},
  series    = {Studies in Applied and Numerical Mathematics},
  year      = {1981},
  address   = {Philadelphia, PA, United States},
  doi       = {10.1137/1.9781611970883}
}

@article{Amestoy_et_al_2001,
  author  = {P. R. Amestoy and I. S. Duff and J. Koster and J.-Y. L'Excellent},
  title   = {A fully asynchronous multifrontal solver using distributed dynamic scheduling},
  journal = {SIAM Journal on Matrix Analysis and Applications},
  volume  = {23},
  number  = {1},
  year    = {2001},
  pages   = {15--41}
}

@phdthesis{Andrews_2025,
  title  = {Geometric numerical integration via auxiliary variables},
  author = {Andrews, B. D.},
  school = {University of Oxford},
  year   = {2025}
}

@misc{Andrews_et_al_2026,
  title  = {Automated {Galerkin} time stepping in {Irksome}},
  author = {Andrews, B. D. and Brubeck, P. and Farrell, P. E. and Kirby, R. C. and MacLachlan, S. P.},
  year   = {2026},
  note   = {arXiv:2606.27300},
  doi    = {10.48550/arXiv.2606.27300}
}

@article{Andrews_Farrell_2025a,
  title   = {Enforcing conservation laws and dissipation inequalities numerically via auxiliary variables},
  volume  = {47},
  issn    = {1064-8275},
  doi     = {10.1137/25M1756673},
  number  = {6},
  journal = {SIAM Journal on Scientific Computing},
  author  = {Andrews, B. D. and Farrell, P. E.},
  month   = dec,
  year    = {2025},
  pages   = {A3516--A3535}
}

@misc{Andrews_Farrell_2026_Zenodo,
  author = {Andrews, B. D. and Farrell, P. E.},
  key    = {zenodo/Zenodo-20260814.1},
  title  = {{Software used in `Conservative and dissipative discretisations of multi-conservative ODEs and GENERIC systems'}},
  month  = aug,
  year   = {2026},
  doi    = {10.5281/zenodo.21932393}
}

@article{Balay_et_al_2024,
  author      = {Balay, S. and Abhyankar, S. and Adams, M. F. and Benson, S. and Brown, J. and Brune, P. and Buschelman, K. and Constantinescu, E. and Dalcin, L. and Dener, A. and Eijkhout, V. and Faibussowitsch, J. and Gropp, W. D. and Hapla, V. and Isaac, T. and Jolivet, P. and Karpeev, D. and Kaushik, D. and Knepley, M. G. and Kong, F. and Kruger, S. and May, D. A. and McInnes, L. C. and Mills, R. T. and Mitchell, L. and Munson, T. and Roman, J. E. and Rupp, K. and Sanan, P. and Sarich, J. and Smith, B. F. and Zampini, S. and Zhang, H. and Zhang, H. and Zhang, J.},
  title       = {{PETSc/TAO} users manual},
  institution = {Argonne National Laboratory},
  journal     = {Argonne National Laboratory},
  number      = {ANL-21/39 - Revision 3.21},
  doi         = {10.2172/2205494},
  year        = {2024}
}

@article{Benjamin_1972,
  title     = {The stability of solitary waves},
  volume    = {328},
  issn      = {0080-4630},
  number    = {1573},
  journal   = {Proceedings of the Royal Society of London. Series A, Mathematical and Physical Sciences},
  publisher = {The Royal Society},
  author    = {Benjamin, T. B.},
  year      = {1972},
  pages     = {153--183}
}

@article{Benjamin_Bona_Mahony_1972,
  title   = {Model equations for long waves in nonlinear dispersive systems},
  volume  = {272},
  doi     = {10.1098/rsta.1972.0032},
  number  = {1220},
  journal = {Philosophical Transactions of the Royal Society of London. Series A, Mathematical and Physical Sciences},
  author  = {Benjamin, T. B. and Bona, J. L. and Mahony, J. J.},
  month   = mar,
  year    = {1972},
  pages   = {47--78}
}

@article{Brugnano_FrascaCaccia_Iavernaro_2019,
  title    = {Line integral solution of {Hamiltonian} {PDEs}},
  volume   = {7},
  issn     = {2227-7390},
  doi      = {10.3390/math7030275},
  language = {en},
  number   = {3},
  journal  = {Mathematics},
  author   = {Brugnano, L. and Frasca-Caccia, G. and Iavernaro, F.},
  month    = mar,
  year     = {2019},
  pages    = {275}
}

@article{Brugnano_Iavernaro_2012,
  title    = {Line integral methods which preserve all invariants of conservative problems},
  volume   = {236},
  issn     = {0377-0427},
  doi      = {10.1016/j.cam.2012.03.026},
  number   = {16},
  journal  = {Journal of Computational and Applied Mathematics},
  author   = {Brugnano, L. and Iavernaro, F.},
  month    = oct,
  year     = {2012},
  pages    = {3905--3919}
}

@book{Brugnano_Iavernaro_2016,
  title     = {Line integral methods for conservative problems},
  isbn      = {978-1-4822-6385-5},
  language  = {en},
  publisher = {CRC Press},
  author    = {Brugnano, L. and Iavernaro, F.},
  month     = mar,
  year      = {2016},
  address   = {Boca Raton, FL, United States},
  doi       = {10.1201/b19319}
}

@article{Cohen_Hairer_2011,
  title    = {Linear energy-preserving integrators for {Poisson} systems},
  volume   = {51},
  issn     = {1572-9125},
  doi      = {10.1007/s10543-011-0310-z},
  language = {en},
  number   = {1},
  journal  = {BIT Numerical Mathematics},
  author   = {Cohen, D. and Hairer, E.},
  month    = mar,
  year     = {2011},
  pages    = {91--101}
}

@article{Dahlby_Owren_Yaguchi_2011,
  title    = {Preserving multiple first integrals by discrete gradients},
  volume   = {44},
  issn     = {1751-8121},
  doi      = {10.1088/1751-8113/44/30/305205},
  language = {en},
  number   = {30},
  journal  = {Journal of Physics A: Mathematical and Theoretical},
  author   = {Dahlby, M. and Owren, B. and Yaguchi, T.},
  month    = jul,
  year     = {2011},
  pages    = {305205}
}

@book{Ern_Guermond_2021a,
  address    = {Cham, Switzerland},
  series     = {Texts in Applied Mathematics},
  title      = {Finite elements {I}: approximation and interpolation},
  volume     = {72},
  isbn       = {978-3-030-56340-0 978-3-030-56341-7},
  shorttitle = {Finite elements {I}},
  language   = {en},
  publisher  = {Springer International Publishing},
  author     = {Ern, A. and Guermond, J.-L.},
  year       = {2021},
  doi        = {10.1007/978-3-030-56341-7}
}

@book{Ern_Guermond_2021c,
  address    = {Cham, Switzerland},
  series     = {Texts in Applied Mathematics},
  title      = {Finite elements {III}: first-order and time-dependent {PDEs}},
  volume     = {74},
  isbn       = {978-3-030-57347-8 978-3-030-57348-5},
  shorttitle = {Finite elements {III}},
  language   = {en},
  publisher  = {Springer International Publishing},
  author     = {Ern, A. and Guermond, J.-L.},
  year       = {2021},
  doi        = {10.1007/978-3-030-57348-5}
}

@article{Evans_1990,
  title     = {Superintegrability in classical mechanics},
  volume    = {41},
  doi       = {10.1103/PhysRevA.41.5666},
  number    = {10},
  urldate   = {2026-07-30},
  journal   = {Physical Review A},
  publisher = {American Physical Society},
  author    = {Evans, N. W.},
  month     = may,
  year      = {1990},
  pages     = {5666--5676}
}

@article{Farrell_Kirby_MarchenaMenendez_2021,
  title   = {{Irksome}: automating {Runge}--{Kutta} time-stepping for finite element methods},
  author  = {Farrell, P. E. and Kirby, R. C. and Marchena-Menendez, J.},
  journal = {ACM Transactions on Mathematical Software},
  volume  = {47},
  number  = {4},
  pages   = {1--26},
  year    = {2021},
  doi     = {10.1145/3466168}
}

@article{Gay-Balmaz_Yoshimura_2017,
  title    = {A {Lagrangian} variational formulation for nonequilibrium thermodynamics. {Part} {I}: {Discrete} systems},
  volume   = {111},
  issn     = {0393-0440},
  doi      = {10.1016/j.geomphys.2016.08.018},
  journal  = {Journal of Geometry and Physics},
  author   = {Gay-Balmaz, F. and Yoshimura, H.},
  month    = jan,
  year     = {2017},
  pages    = {169--193}
}

@article{Giesselmann_Karsai_Tscherpel_2025,
  title   = {Energy-consistent {Petrov}--{Galerkin} time discretization of port-{Hamiltonian} systems},
  volume  = {11},
  issn    = {2426-8399},
  doi     = {10.5802/smai-jcm.127},
  journal = {SMAI Journal of Computational Mathematics},
  author  = {Giesselmann, J. and Karsai, A. and Tscherpel, T.},
  year    = {2025},
  pages   = {335--367}
}

@article{Gonzalez_1996,
  title    = {Time integration and discrete {H}amiltonian systems},
  volume   = {6},
  issn     = {1432-1467},
  doi      = {10.1007/BF02440162},
  language = {en},
  number   = {5},
  journal  = {Journal of Nonlinear Science},
  author   = {Gonzalez, O.},
  month    = sep,
  year     = {1996},
  pages    = {449--467}
}

@article{Grillakis_Shatah_Strauss_1987,
  title   = {Stability theory of solitary waves in the presence of symmetry, {I}},
  volume  = {74},
  issn    = {0022-1236},
  doi     = {10.1016/0022-1236(87)90044-9},
  number  = {1},
  urldate = {2026-08-05},
  journal = {Journal of Functional Analysis},
  author  = {Grillakis, M. and Shatah, J. and Strauss, W.},
  month   = sep,
  year    = {1987},
  pages   = {160--197}
}

@article{Grmela_1984,
  title   = {Particle and bracket formulations of kinetic equations},
  volume  = {28},
  journal = {Contemporary Mathematics},
  author  = {Grmela, M.},
  year    = {1984},
  pages   = {125--132}
}

@article{Grmela_Ottinger_1997,
  title    = {Dynamics and thermodynamics of complex fluids. {I}. {Development} of a general formalism},
  volume   = {56},
  doi      = {10.1103/PhysRevE.56.6620},
  number   = {6},
  journal  = {Physical Review E},
  author   = {Grmela, M. and Öttinger, H. C.},
  month    = dec,
  year     = {1997},
  pages    = {6620--6632}
}

@article{Hairer_et_al_2006,
  title    = {Geometric Numerical Integration},
  volume   = {3},
  issn     = {1660-8933},
  doi      = {10.4171/owr/2006/14},
  language = {en},
  number   = {1},
  journal  = {Oberwolfach Reports},
  author   = {Hairer, E. and Hochbruck, M. and Iserles, A. and Lubich, C.},
  month    = dec,
  year     = {2006},
  pages    = {805--882}
}

@book{Hairer_Lubich_Wanner_2006,
  title      = {Geometric numerical integration: structure-preserving algorithms for ordinary differential equations},
  isbn       = {978-3-540-30666-5},
  shorttitle = {Geometric Numerical Integration},
  language   = {en},
  publisher  = {Springer Science \& Business Media},
  author     = {Hairer, E. and Lubich, C. and Wanner, G.},
  month      = may,
  year       = {2006},
  address    = {Heidelberg, Germany},
  doi        = {10.1007/3-540-30666-8}
}

@article{Ham_et_al_2023,
  title        = {Firedrake user manual},
  author       = {Ham, D. A. and Kelly, P. H. J. and Mitchell, L. and Cotter, C. J. and Kirby, R. C. and Sagiyama, K. and Bouziani, N. and Vorderwuelbecke, S. and Gregory, T. J. and Betteridge, J. and Shapero, D. R. and Nixon-Hill, R. W. and Ward, C. J. and Farrell, P. E. and Brubeck, P. D. and Marsden, I. and Gibson, T. H. and Homolya, M. and Sun, T. and McRae, A. T. T. and Luporini, F. and Gregory, A. and Lange, M. and Funke, S. W. and Rathgeber, F. and Bercea, G.-T. and Markall, G. R.},
  organization = {Imperial College London and University of Oxford and Baylor University and University of Washington},
  journal      = {Imperial College London, University of Oxford, Baylor University, University of Washington},
  year         = {2023},
  doi          = {10.25561/104839}
}

@article{Hulme_1972a,
  title     = {One-step piecewise polynomial {Galerkin} methods for initial value problems},
  volume    = {26},
  issn      = {0025-5718},
  doi       = {10.2307/2005168},
  number    = {118},
  journal   = {Mathematics of Computation},
  publisher = {American Mathematical Society},
  author    = {Hulme, B. L.},
  year      = {1972},
  pages     = {415--426}
}

@article{Hulme_1972b,
  title   = {Discrete {Galerkin} and related one-step methods for ordinary differential equations},
  volume  = {26},
  issn    = {0025-5718},
  doi     = {10.1090/S0025-5718-1972-0315899-8},
  number  = {120},
  journal = {Mathematics of Computation},
  author  = {Hulme, B. L.},
  year    = {1972},
  pages   = {881--891}
}

@article{Kaufman_1984,
  title      = {Dissipative {Hamiltonian} systems: {A} unifying principle},
  volume     = {100},
  issn       = {0375-9601},
  shorttitle = {Dissipative {Hamiltonian} systems},
  doi        = {10.1016/0375-9601(84)90634-0},
  number     = {8},
  journal    = {Physics Letters A},
  author     = {Kaufman, A. N.},
  month      = feb,
  year       = {1984},
  pages      = {419--422}
}

@article{Kovalevskaya_1889,
  title    = {Sur le problème de la rotation d'un corps solide autour d'un point fixe},
  volume   = {12},
  language = {fr},
  journal  = {Acta Mathematica},
  author   = {Kovalevskaya, S.},
  year     = {1889},
  pages    = {177--232},
  doi      = {10.1007/BF02592182}
}

@article{LaBudde_Greenspan_1974,
  title    = {Discrete mechanics---a general treatment},
  volume   = {15},
  issn     = {0021-9991},
  doi      = {10.1016/0021-9991(74)90081-3},
  number   = {2},
  journal  = {Journal of Computational Physics},
  author   = {LaBudde, R. A. and Greenspan, D.},
  month    = jun,
  year     = {1974},
  pages    = {134--167}
}

@misc{Lombardi_Pagliantini_2025,
  title     = {Conformal variational discretisation of infinite dimensional {Hamiltonian} systems with gradient flow dissipation},
  doi       = {10.48550/arXiv.2412.06310},
  publisher = {arXiv},
  author    = {Lombardi, D. and Pagliantini, C.},
  month     = may,
  year      = {2025}
}

@article{McLachlan_Quispel_Robidoux_1999,
  title    = {Geometric integration using discrete gradients},
  volume   = {357},
  doi      = {10.1098/rsta.1999.0363},
  number   = {1754},
  journal  = {Philosophical Transactions of the Royal Society of London. Series A: Mathematical, Physical and Engineering Sciences},
  author   = {McLachlan, R. I. and Quispel, G. R. W. and Robidoux, N.},
  month    = apr,
  year     = {1999},
  pages    = {1021--1045}
}

@article{Miller_Post_Winternitz_2013,
  title     = {Classical and quantum superintegrability with applications},
  volume    = {46},
  issn      = {1751-8121},
  doi       = {10.1088/1751-8113/46/42/423001},
  language  = {en},
  number    = {42},
  journal   = {Journal of Physics A: Mathematical and Theoretical},
  publisher = {IOP Publishing},
  author    = {Miller, W. and Post, S. and Winternitz, P.},
  month     = oct,
  year      = {2013},
  pages     = {423001}
}

@article{Morrison_1984a,
  title   = {Bracket formulation for irreversible classical fields},
  volume  = {100},
  issn    = {0375-9601},
  doi     = {10.1016/0375-9601(84)90635-2},
  number  = {8},
  journal = {Physics Letters A},
  author  = {Morrison, P. J.},
  month   = feb,
  year    = {1984},
  pages   = {423--427}
}

@techreport{Morrison_1984b,
  title       = {Some observations regarding brackets and dissipation},
  number      = {PAM-228},
  institution = {University of California at Berkeley},
  author      = {Morrison, P. J.},
  year        = {1984}
}

@article{Morrison_1986,
  title   = {A paradigm for joined {Hamiltonian} and dissipative systems},
  volume  = {18},
  issn    = {0167-2789},
  doi     = {10.1016/0167-2789(86)90209-5},
  number  = {1},
  journal = {Physica D: Nonlinear Phenomena},
  author  = {Morrison, P. J.},
  year    = {1986},
  pages   = {410--419}
}

@article{Nambu_1973,
  title   = {Generalized {Hamiltonian} dynamics},
  volume  = {7},
  doi     = {10.1103/PhysRevD.7.2405},
  number  = {8},
  journal = {Physical Review D},
  author  = {Nambu, Y.},
  year    = {1973},
  pages   = {2405--2412}
}

@article{Ottinger_1997,
  title    = {{GENERIC} formulation of {Boltzmann's} kinetic equation},
  volume   = {22},
  issn     = {1437-4358},
  doi      = {10.1515/jnet.1997.22.4.386},
  journal  = {Journal of Non-Equilibrium Thermodynamics},
  number   = {4},
  author   = {Öttinger, H. C.},
  month    = jan,
  year     = {1997},
  pages    = {386--391}
}

@book{Ottinger_2005,
  title     = {Beyond {Equilibrium} {Thermodynamics}},
  isbn      = {978-0-471-72791-0},
  publisher = {John Wiley \& Sons},
  author    = {Öttinger, H. C.},
  month     = may,
  year      = {2005},
  address   = {New York, NY, United States},
  doi       = {10.1002/0471727903}
}

@article{Ottinger_Grmela_1997,
  title    = {Dynamics and thermodynamics of complex fluids. {II}. {Illustrations} of a general formalism},
  volume   = {56},
  doi      = {10.1103/PhysRevE.56.6633},
  number   = {6},
  journal  = {Physical Review E},
  author   = {Öttinger, H. C. and Grmela, M.},
  month    = dec,
  year     = {1997},
  pages    = {6633--6655}
}

@article{Romero_2009,
  title    = {Thermodynamically consistent time-stepping algorithms for non-linear thermomechanical systems},
  volume   = {79},
  doi      = {10.1002/nme.2588},
  number   = {6},
  journal  = {International Journal for Numerical Methods in Engineering},
  author   = {Romero, I.},
  month    = mar,
  year     = {2009},
  pages    = {706--732}
}

@article{Schulz_Wan_2025,
  title   = {Minimal $\ell^2$ norm discrete multiplier method},
  volume  = {12},
  doi     = {10.3934/jcd.2024029},
  number  = {2},
  journal = {Journal of Computational Dynamics},
  author  = {Schulz, E. and Wan, A. T. S.},
  year    = {2025},
  pages   = {212--238}
}

@article{Shang_Ottinger_2020,
  title   = {Structure-preserving integrators for dissipative systems based on reversible--irreversible splitting},
  volume  = {476},
  doi     = {10.1098/rspa.2019.0446},
  number  = {2234},
  journal = {Proceedings of the Royal Society A},
  author  = {Shang, X. and Öttinger, H. C.},
  month   = feb,
  year    = {2020}
}

@book{Taff_1985,
  title     = {Celestial mechanics},
  isbn      = {978-0-471-89316-5},
  language  = {en},
  publisher = {Wiley},
  author    = {Taff, L. G.},
  month     = may,
  year      = {1985},
  address   = {New York, NY, United States}
}

@article{Takhtajan_1994,
  title   = {On foundation of the generalized {Nambu} mechanics},
  volume  = {160},
  issn    = {1432-0916},
  doi     = {10.1007/BF02103278},
  number  = {2},
  journal = {Communications in Mathematical Physics},
  author  = {Takhtajan, L.},
  year    = {1994},
  pages   = {295--315}
}

@book{Tu_2010,
  title     = {An introduction to manifolds},
  isbn      = {978-1-4419-7400-6},
  language  = {en},
  publisher = {Springer Science \& Business Media},
  author    = {Tu, L. W.},
  month     = oct,
  year      = {2010},
  address   = {New York, NY, United States},
  doi       = {10.1007/978-1-4419-7400-6}
}

@article{Wan_Bihlo_Nave_2017,
  title   = {Conservative methods for dynamical systems},
  volume  = {55},
  doi     = {10.1137/16M110719X},
  number  = {5},
  journal = {SIAM Journal on Numerical Analysis},
  author  = {Wan, A. T. S. and Bihlo, A. and Nave, J.-C.},
  year    = {2017},
  pages   = {2255--2285}
}

@article{Wan_Nave_2018,
  title   = {On the arbitrarily long-term stability of conservative methods},
  volume  = {56},
  doi     = {10.1137/16M1085929},
  number  = {5},
  journal = {SIAM Journal on Numerical Analysis},
  author  = {Wan, A. T. S. and Nave, J.-C.},
  year    = {2018},
  pages   = {2751--2775}
}

@article{Zhong_Marsden_1988,
  title    = {{Lie--Poisson Hamilton--Jacobi theory and Lie--Poisson integrators}},
  volume   = {133},
  issn     = {0375-9601},
  doi      = {10.1016/0375-9601(88)90773-6},
  number   = {3},
  journal  = {Physics Letters A},
  author   = {Ge, Z. and Marsden, J. E.},
  month    = nov,
  year     = {1988},
  pages    = {134--139}
}

\end{document}